\documentclass[english]{article}
\usepackage[T1]{fontenc}
\usepackage{algorithmic}
\usepackage{geometry}
\usepackage{float}
\usepackage{amsthm}
\usepackage{mathtools}
\usepackage{amsmath}
\usepackage{amssymb}
\usepackage{setspace}
\usepackage{bbm}
\usepackage{esint}
\usepackage{comment}
\usepackage{cleveref}
\usepackage{enumitem}
\usepackage{graphicx}
\usepackage{caption}
\usepackage{thmtools}
\usepackage{thm-restate}
\usepackage{subcaption}

\makeatletter

\newcommand{\bbP}{\mathbb{P}}
\newcommand{\bbR}{\mathbb{R}}
\newcommand{\bbN}{\mathbb{N}}
\newcommand{\bbE}{\mathbb{E}}
\newcommand{\bbB}{\mathbb{B}}
\newcommand{\bbone}{\mathbbm{1}}
\newcommand{\calH}{\mathcal{H}}
\newcommand{\underlambda}{\underline{\lambda}}
\newcommand{\overlambda}{\overline{\lambda}}
\newcommand{\overpi}{\overline{\pi}}
\newcommand{\underpi}{\underline{\pi}}

\newcommand{\per}{\mathrm{Per}}
\floatstyle{ruled}
\newfloat{algorithm}{tbp}{loa}
\providecommand{\algorithmname}{Algorithm}
\floatname{algorithm}{\protect\algorithmname}

\numberwithin{equation}{section}
\numberwithin{figure}{section}
\theoremstyle{plain}
\newtheorem{thm}{\protect\theoremname}
  \theoremstyle{plain}
  \newtheorem{lem}[thm]{\protect\lemmaname}

  \theoremstyle{plain}

  \theoremstyle{plain}

\makeatother

\usepackage{babel}
  \providecommand{\lemmaname}{Lemma}
  \providecommand{\corollaryname}{Corollary}
  \providecommand{\propname}{Proposition}
\providecommand{\theoremname}{Theorem}

\begin{document}
\begin{doublespace}
\end{doublespace}

\title{Multi-Step Bayesian Optimization for One-Dimensional Feasibility Determination}
\author{J. Massey Cashore \and Lemuel Kumarga \and Peter I. Frazier}
\date{\vspace{-5ex}}
\maketitle
\begin{abstract}
Bayesian optimization methods allocate limited sampling budgets to maximize expensive-to-evaluate functions.  One-step-lookahead policies are often used, but computing optimal multi-step-lookahead policies remains a challenge.  We consider a specialized Bayesian optimization problem: finding the superlevel set of an expensive one-dimensional function, with a Markov process prior.  We compute the Bayes-optimal sampling policy efficiently, and characterize the suboptimality of one-step lookahead.
Our numerical experiments demonstrate that the one-step lookahead policy is close to optimal in 
this problem, performing within 98\% of optimal in the experimental settings considered.
\end{abstract}

\newcommand{\pfcomment}[1]{ PF: {\it #1}}
\newcommand{\mccomment}[1]{ MC: {\it #1}}
\section{Introduction}
\label{introsection}

We consider the problem of adaptively allocating sampling effort to efficiently estimate sub- and super-level sets of a one-dimensional Markov process, or more general additive functionals of this process.  We use a decomposition property to show how the optimal procedure may be computed efficiently, circumventing the curse of dimensionality.  We then use our ability to compute the optimal policy to study the suboptimality gap of commonly used one-step lookahead procedures in this problem.

The problem we consider falls within the class of problems considered in the large and rapidly growing literature on Bayesian optimization \cite{Ku64,Mo89,JoScWe98,FoSoKe08,snoek2012practical}, which seeks to develop adaptive sampling algorithms that estimate functionals, especially the location of a global maximum, of some underlying and unknown function in a query efficient way.  Such problems arise when optimizing an objective that is computed via a long-running computer code \cite{FoSoKe08,snoek2012practical} or some other expensive process \cite{BrCoFr09,FrazierWang2014} that severely limits the number of times it may be sampled.  This literature places a Bayesian prior distribution on the underlying function, and views it as a realization of a stochastic process, most frequently a Gaussian process.

In such problems, a Bayes-optimal algorithm is one that minimizes the expected loss under the prior suffered from mis-estimation of the underlying functional of interest, where the cost of sampling is either factored directly into the objective (as considered by \cite{ChickFrazier2012,Betro1991}), or a sampling budget is enforced as a constraint (as considered by \cite{GinsbourgerLeRiche2010,Frazier2012}). When optimization is the goal, this loss function is the opportunity cost --- the difference in value between the point that is believed to be the best, and the value of the true global optimum --- but when other functionals are of interest another loss function may be appropriate.

In principle, a Bayes-optimal algorithm may be computed using stochastic dynamic programming by understanding that this problem is a partially observable Markov decision process (POMDP) \cite{GinsbourgerLeRiche2010}.  However, the curse of dimensionality \cite{Po07} prevents actually computing the solution through brute-force approaches.

Thus, almost all of the literature has focused on approximate schemes, which in many cases are inspired by this view of the problem as a partially observable Markov decision process, but that do not actually solve the POMDP.  Two commonly used methods of this type are the expected improvement method \cite{Mo89,JoScWe98} and the knowledge-gradient method \cite{FrazierPowellDayanik2009,ScottFrazierPowell2011}, which use one-step lookahead approaches, based on different assumptions about what points are eligible for selection once sampling stops \cite{FrazierWang2014}.  Two-step lookahead approaches have also been implemented computationally in \cite{BetroSchoen1991,GinsbourgerLeRiche2010}.

In contrast, we focus on calculating the Bayes-optimal algorithm. Our primary contribution is to  show that it can be computed efficiently in Bayesian optimization problems that satisfy four assumptions:
\begin{itemize}
    \item the underlying function is one-dimensional, as considered by \cite{Ku64,CaZi99,CaZi00,Calvin2001,CaZi02,Perttunen1991,Locatelli1997}.
    \item the Bayesian prior on this function has the Markov property (e.g., a Wiener process prior, as used by \cite{Ku64,Perttunen1991,Archetti1979,Zilinskas1985,Ritter1990,LocatelliSchoen1995,Locatelli1997,Calvin2001}, or an Ornstein-Uhlenbeck prior \cite{PerttunenStuckman1989,PerttunenStuckman1990}). 
\item the loss function is additive across location, as arises when the goal is to determine feasibility of points, as in \cite{Gardner2014}, or to determine the set of points that are better than some known standard, as in \cite{XieFrazier2013a}.
\item the limit on sampling is imposed as an additive cost in the objective, as in \cite{ChickFrazier2012,Betro1991}, or as a constraint on the expected number of samples taken, as in \cite{GinsbourgerLeRiche2010,Frazier2012}.
\end{itemize}

As a second contribution, we also provide an upper bound on the value of the Bayes optimal policy when the limit on sampling is imposed as an almost sure constraint on the number of samples taken.

While one dimensional feasibility determination problems do arise in practice \cite{hitz14fully}, and we expect that the optimal policy can provide a great deal of value in those settings, a large fraction of practical Bayesian optimization problems violate one or more of the assumptions above, because many problems are in more than one dimension, and because non-Markov Gaussian processes are often used as priors \cite{gotovos13active, Sa02,BrCoFr09}.  Optimization is also a more common goal in the literature than super-level set determination (though in practice it is often just as useful to provide a set of points that perform well, i.e., that reside in some super-level set, from which a final decision can be selected based on other criteria).

Thus, we view our primary contribution as providing a specialized but nevertheless rich class of Bayesian optimization problems on which the performance of widely applicable heuristic procedures, such as the one-step lookahead procedures described above, may be studied relative to Bayes-optimal procedures.  This guides algorithm development --- if a heuristic procedure performs close to optimal on a set of problems, then this suggests that further improvement is not necessary even for other similar problems for which the optimality gap cannot be evaluated.  In contrast, if all known heuristic procedures perform substantially worse than optimal on a set of problems, then this suggests that further algorithm development is worthwhile.

There is some complementary theoretical analysis in the literature of Bayes-optimal procedures for this and related problems.  
Much of it focuses on asymptotic analyses, and includes proofs of consistency for the Efficient Global Optimization (EGO) \cite{VazquezBect2010} and P algorithms \cite{CaZi99}, 
as well as convergence rates for these and closely related algorithms \cite{CaZi00,Bull2011,Calvin2001}.
In terms of finite-time analyses, \cite{Srinivas2010,Grunewalder2010} provide regret bounds for the closely related problem of Bayesian optimization in the bandit setting,
but  while these bounds characterize performance, slack in the bounds' constants creates a potentially large multiplicative gap in which performance may lie.
In the problems of multiple comparisons with a known standard and stochastic root-finding, procedures for computing explicit Bayes optimal procedures have been developed \cite{XieFrazier2013a,WaeberFrazierHenderson2013}, but these problems are only distantly related to Bayesian optimization.
Thus, exact performance of optimal finite-time procedures has remained unknown in Bayesian optimization.

Below, in Section \ref{problemdescsec}, we provide a formal description of the problem.  Our main
 results are in Section \ref{mainresultssec}, where we significantly reduce the state-space for
a dynamic program giving rise to a Bayes-optimal policy.
In Section \ref{upperboundsec} we consider the relationship between the cost-per-sample setting
and the constrained-budget setting, showing how the optimal value for the former can be used
to compute the optimal value of the latter.
In Section \ref{simulationsec} we present numerical results, illustrating the behavior of the 
optimal policy and using it to analyze the optimality gap for a one-step lookahead procedure.  
Finally, in Section \ref{conclusionsec}, we conclude.

\newcommand{\zap}[1]{}

\section{Problem Description}
\label{problemdescsec}

Let $Y = (Y(x) : x\ge 0)$ be a Markov process over the positive real line, and let $[a,b]$ be a given interval, $0<a<b<\infty$.  We consider adaptive sampling policies that characterize $Y$ over $[a,b]$.

We will consider histories of the form $\{ (x_t, Y(x_t)) : t=1,\ldots,T \}$ for some sequence of adaptively chosen points $(x_t : t=1,\ldots, T)$ at which measurements occur.

Let $\calH = \cup_{T=0}^\infty (\bbR_+ \times \bbR)^T$ be the space of all possible histories. A policy $\pi : \calH \mapsto \bbR_+ \cup \{\Delta\}$ is a measurable function that maps the current history to either a point to be sampled next, or to the symbol $\Delta$, which indicates the decision to stop.
We let $\Pi$ indicate the space of all such policies.

We begin with an initial history $H_0\in\calH$. For simplicity of analysis we assume 
that $H_0$ contains endpoint observations, that is $(a,Y(a)), (b,Y(b))\in H_0$, but 
our results can be extended to the case where it does not.
We define histories $H_t$ and 
decisions $x_t$ recursively, letting $x_{t+1}=\pi(H_t)$ and letting 
\begin{equation}
\label{updateequation}
H_{t+1} \coloneqq \begin{cases}
H_t \cup \{(x_{t+1},Y(x_{t+1}))\},&\text{if $x_{t+1}\ne\Delta$,}\\
H_t, &\text{if $x_{t+1}=\Delta$,}
\end{cases}
\end{equation}
so that the point sampled and the resulting observation of $Y$ is added to the history if the policy chooses to sample,
and the history remains unchanged once the policy chooses to stop sampling.
As indicated, if a policy measures a point already in the history, then the history remains unchanged.
Also note that because the history is a set of tuples, the policy's next action cannot depend on the order
in which observations were made.

We define
$\tau = \inf\left\{ t\ge0:\pi(H_t)=\triangle\right\}$
to be the total number of samples taken by a policy.
When necessary, we will write $\tau^\pi$ to emphasize
the policy on which $\tau$ depends. 

We also define $\bbP^\pi$ to be the distribution over histories with respect
to the randomness in $Y$ and the decisions made by $\pi$, for any $\pi\in\Pi$.
We let $\bbE^\pi$ denote the expectation with respect to this distribution.

We seek to characterize $Y$ by assigning each point $x\in [a,b]$ a label, or class,
based on our knowledge of $Y(x)$.
Suppose there are $n<\infty$ classes to which a point might
belong and let $I$ be an index set such that each element corresponds to
one class.  At time $\tau$, we will use the information collected, encoded in $H_\tau$, 
to classify each point in $[a,b]$: based on $H_\tau$ we construct a partition $\{B_i : i \in I\}$
of $[a,b]$, such that each $B_i$ is a measurable subset of $[a,b]$.   
If $x\in B_i$, we say that $x$ belongs to the $i$th class.
We will receive a reward $R_{[a,b]}(H_\tau)$, defined below, that depends on
the accuracy of this classification.

To formalize this we first choose bounded measurable functions 
$f_i : \bbR \to \bbR$ for each $i\in I$.
The function $f_i$ is meant to reward or penalize the classification $x\in B_i$
given the true value $Y(x)$.
In addition to requiring the $f_i$ be measurable and bounded, 
we also require that, for each $i\in I$,
$f_i$ satisfies the following inequality: 
\begin{equation}
\label{rewardfunctioninequality}
	\int_{[a,b]} \bbE \left[ \left|f_i\circ Y(x)\right| \right]dx < \infty.
\end{equation}
By Fubini's theorem, this inequality will allow us below to interchange the integral over the domain
of $Y$ and the expectation over the randomness in $Y$.

Now, fix any partition $\bbB=\{B_i : i\in I\}$ and $H\in\calH$. We define the expected
reward of the partition $\bbB$ given $H_\tau = H$ over $[a,b]$ to be
\begin{equation}
R_{[a,b]}(H, \bbB) = \bbE\left[ \sum_{i\in I} \int_{B_i} f_i \circ Y(x)dx \mid H\right]
 =  \sum_{i\in I} \int_{B_i} \bbE \left[f_i \circ Y(x)\mid H \right]dx,
\label{nicerewardexp}
\end{equation}
where the last equality holds due to (\ref{rewardfunctioninequality}) and Fubini's theorem.
Observing (\ref{nicerewardexp}), a partition maximizing $R_{[a,b]}(H, \bbB)$
is any $\bbB^* = \{B^*_i : i \in I\}$ such that for all $j\in I$, if $x\in B_j$, then
\[
j \in \underset{i}{\operatorname{argmax}}\mbox{ }\bbE\left[f_i \circ Y(x)\mid H \right].
\]
That is, $x$ belongs to any class $j$ maximizing $\bbE\left[f_j \circ Y(x)\mid H \right]$.  

We define the expected reward $R_{[a,b]}(H)$ for any $H\in\calH$ to be the expected reward
of any optimal partition given $H$. That is, 
\begin{equation}
R_{[a,b]}(H)=R_{[a,b]}(H,\bbB^*)=
\int_{[a,b]}\max_i\bbE\left[f_i \circ Y(x) \mid H \right]dx.
\end{equation}
Note that because we choose the functions $\{f_i:i\in I\}$ to be bounded, it follows
that there exists some constant $C$ such that 
\begin{equation}
\label{rewardupperbound}
|R_{[a,b]}(H)|\le C (b-a),
\end{equation}
for every $H\in \calH$.

We now describe how this framework can be specialized to the problem of 
estimating superlevel sets.
Recall the superlevel set of a function $g:[a,b]\to\bbR$ with respect to
the threshold $k$ is the set $\{x\in[a,b]:g(x)\ge k\}$. 
In this context we use the index set $I=\{+,-\}$
corresponding to the classification of a point as above or below the threshold.
We give two reasonable choices for the functions $f_+$ and $f_-$:
\begin{enumerate}
\item $f_+(y)=\bbone\{y\ge k\}$ and $f_-(y)=\bbone\{y\le k\}$. 
These functions are clearly measurable, bounded, and satisfy the inequality (\ref{rewardfunctioninequality}).
\item 
\begin{align*}
f_+(y) &= 
\begin{cases}
y-k&\text{if $|y-k|\le C$,}\\
C\cdot\mathrm{sign}(y-k) &\text{otherwise,}
\end{cases}
&& \text{and} &
f_-(y) &= 
\begin{cases}
k-y&\text{if $|k-y|\le C$,}\\
C\cdot\mathrm{sign}(k-y) &\text{otherwise,}
\end{cases}
\end{align*}
for some constant $C>0$ chosen a priori.
We only consider these reward functions for Markov processes $Y$ such that
the inequality (\ref{rewardfunctioninequality}) is satisfied.
\end{enumerate}

Although we focus on superlevel set detection, this framework can
be used to classify points based on other properties of $Y(x)$.  For example,
we could consider two thresholds for the range of $Y(x)$, and classify each point
as being below both, above both, or in between them.

The performance of a policy $\pi$
at state $H$ over $[a,b]$ is the expected value of the final reward
less the cost associated with the expected number of samples starting
from an initial history $H$. That is, given a cost-per-sample of $c>0$,
the performance is defined as:
\begin{equation}
\label{performanceeqfirst}
\mathrm{Per}(\pi,c,H) = \bbE^\pi \left[ R_{[a,b]}(H_\tau) - c\tau \mid H\right].
\end{equation}
As a consequence of (\ref{rewardupperbound}), a policy that has non-zero probability of
taking infinitely many samples at a state $H$ (i.e. $\bbP^\pi(\tau=\infty|H)\ne 0$) 
achieves a performance of $-\infty$ at $H$.

Finally, the value of a state $H$ is defined to be the supremum of the performance 
over all policies. That is,
\begin{equation}
\label{valuefunctiondef}
V_{[a,b]}(H) = \sup_{\pi \in \Pi} 
\bbE^\pi \left[ R_{[a,b]}(H_\tau) - c\tau \mid H\right].
\end{equation}
We call (\ref{valuefunctiondef}) the
cost-per-sample setting.  Below, in section \ref{upperboundsec}, we consider
two other related settings: budget-constrained and expected-budget-constrained.

Now, in Section \ref{mainresultssec}, we show how to approximately compute
$\epsilon$-optimal policies for the
cost-per-sample setting when summability is satisfied, i.e., one that attains within $\epsilon$ the supremum in (\ref{valuefunctiondef}) for
the initial history $H_0$ and any $\epsilon>0$. We refer to such an optimal policy as "$\epsilon$-Bayes-optimal" because it is
$\epsilon$-optimal with respect to an expectation taken over the probability distribution of $Y$, which can
be understood to be a Bayesian prior distribution.

\section{Main Results for Cost-Per-Sample Case}
\label{mainresultssec}
To compute a Bayes-optimal policy, we focus on
efficiently computing the value function defined in (\ref{valuefunctiondef}).
Naive dynamic programming can be used, but the dimensionality of the portion
of the state space reachable after $t$ samples grows linearly in $t$,
causing the volume of the state space, and thus the memory and computation required for
dynamic programming, to grow exponentially.  Our main result decomposes the value function,
showing that it is completely determined by its value on some $4$-dimensional set, leading
to its computation as a tractable dynamic program over a state space of small constant dimension. 
In particular, we prove the following:

\begin{restatable}{thm}{decomp}
\label{decompositiontheorem}
Fix interval $[a,b]$ and let $H\in\calH$ be such that observations of $a$ and $b$ are included. 
Let $x_1,\dots,x_{t+2}$ be the observations in $H$ contained within $[a,b]$. Suppose they are
ordered such that $x_i<x_{i+1}$ for all $1\le i<t+2$, so $x_1=a$ and $x_{t+2}=b$.
Define $H_i = \{(x_i, y_i), (x_{i+1}, y_{i+1})\}$
for each $1\le i < t+2$. Then
\begin{equation}
\label{valuedecomp}
V_{[a,b]}(H) = \sum_{i=1}^{t+1} V_{[x_i, x_{i+1}]}(H_i).
\end{equation}  
\end{restatable}

The importance of this theorem is that $V_{[a,b]}(H)$ is completely determined
by its values on $\{H\in\calH : |H|=2\}$, greatly reducing the effective dimension of the dynamic
program's state space.  We show below how this dimension reduction can be used in a recursive
algorithm over a $4$-dimensional state space (the set of histories of length two) to find
$V_{[a,b]}(H)$. Recall that knowledge of the value function $V_{[a,b]}$ at every state
can lead to $\epsilon$-optimal policies. Indeed, if $\epsilon=\kappa_1 + \kappa_2 + \dots$,
and $\pi$ is a policy that, at the $t$th step, selects a point $x_t$ to sample that 
is within $\kappa_t$ of the optimal, then $\pi$ is $\epsilon$-optimal. The point $x_t$ can be
any point such that $V_{[a,b]}(H_t) -\kappa_t \le \bbE[V_{[a,b]}(H_t\cup\{(x_t,Y(x_t))\})|H_t]$.
For further details, see \cite{DynkinYushkevich} section 5.

We can further reduce the state space when $Y$ satisfies additional structure. 
For any $\ell\in\bbR$, define the shift operator $T_\ell:\bbR\to\bbR$ 
by $T_\ell(x)=x+\ell$. We will apply $T_\ell$ to 
elements of $\calH$, and adopt the convention that
$T_\ell(H)=\{(x+\ell,y):(x,y)\in H\}$, i.e. $T_\ell$ only translates
the location of the observations in $H$, and not their values.
We say the Markov process $Y$ is \emph{translation invariant}
if, for any $H\in\calH$, $y\in\bbR$, $x\in\bbR_+$ and $\ell\in\bbR$ such that $x+\ell\ge 0$,
\begin{equation}
\label{translationinvariance}
\bbP(Y(x)\in dy\mid H) = \bbP(Y(x+\ell)\in dy\mid T_\ell(H)).
\end{equation}
The following proposition establishes that if the Markov process is translation
invariant, so is the value function.

\begin{restatable}{prop}{translation}
\label{translationthm}
Suppose $Y$ is translation invariant. Fix interval $[a,b]$ and pick any $\ell\in\bbR$
 such that $a+\ell\ge 0$. Pick any history $H\in\calH$ and let $H'=T_{\ell}(H)$. 
Then
\begin{equation}
\label{translationinv}
V_{[a,b]}(H) = V_{[a',b']}(H')
\end{equation}
where $a'=a+\ell$ and $b'=b+\ell$.
\end{restatable}

Thus when $Y$ satisfies translation invariance, the value function is completely
determined by its values on $\{H\in\calH : |H|=2, (0,y_0)\in H\}$. (The  choice of
$0$ in the $(0,y_0)\in H$ condition is arbitrary; one may replace $0$ by any
other constant in the domain of $Y$). In this case,
$V_{[a,b]}$ can be computed as the result of a dynamic-programming-like recursion
over a $3$-dimensional, rather than $4$-dimensional, state space, as described below. This
reduction in dimension enables faster computation with less memory.

To prove our main results, we first state two technical lemmas:

\begin{lem}
\label{valueoutsideinterval}
Let $H\in \calH$ contain t observations. Let
$H^I=\{(x,y)\in H : x\in [a,b]\}$ denote the set of initial observations
inside $[a,b]$. Then
\begin{equation}
V_{[a, b]}(H) = V_{[a, b]}(H^I).
\end{equation}
\end{lem}

For the following lemma and rest of this section we adopt the notation that,
for $H\in\calH$ and $A\subseteq\bbR$, $H\cap A=\{(x,y)\in H:x\in A\}$.
We will also write $x\in H$ to mean there exists $y\in\bbR$ such that
$(x,y)\in H$.

\begin{lem}
\label{nicepolicieslemma}
Fix interval $[a,b]$.
Let
\begin{itemize}
\item $\Pi^1=\{\pi\in\Pi : \bbP^\pi(\tau<\infty\mid H) = 1 \mbox{ } \forall H\in\calH\}$
 be the set of policies
that almost surely take finitely many samples.
\item $\Pi^2_{[a,b]}=\{\pi\in\Pi : \pi(H)\in [a,b] \mbox{ }\forall H\in\calH\}$ be the set
of policies that only take samples in $[a,b]$.
\item $\Pi^3 = \{\pi\in\Pi : \pi(H)\ne x \mbox{ if } x\in H\mbox{ }\forall H\in\calH\}$
be the set of policies that do not sample the same point twice.
\end{itemize}
Let $\bar{\Pi}_{[a,b]}=\Pi^1 \cap \Pi^2_{[a,b]} \cap \Pi^3$. For all $H\in\calH$, define
\begin{equation}
\bar{V}_{[a,b]}(H)=\sup_{\pi\in\bar{\Pi}_{[a,b]}}\bbE^\pi\left[R_{[a,b]}(H_\tau)-c\tau\mid H\right].
\end{equation}
Then $V_{[a,b]}(H)=\bar{V}_{[a,b]}(H)$ for all $H\in\calH$.
\end{lem}

Lemma \ref{valueoutsideinterval} says that if $H\in\calH$ contains endpoint observations then
the only points in $H$ that affect $V_{[a,b]}(H)$
are those within $[a,b]$.
Lemma \ref{nicepolicieslemma} constructs a subset of $\Pi$ containing an optimal
policy.
The proofs of the above lemmas, as well as Proposition \ref{translationthm} are contained in the appendix. 
We are now in a position to prove the main decomposition theorem.

\begin{proof}[\textbf{Proof of Theorem \ref{decompositiontheorem}}]

We proceed by induction on $t$. When $t=0$ the summation contains
only one term and the result is established. 
Fix $t>0$ and suppose the decomposition (\ref{valuedecomp}) holds
for any $|H|<t+2$.

Define $\tau_{A}$ with respect to any policy $\pi$ to be the number of
points $\pi$ chooses to sample inside the set $A$, for some $A\subseteq [a,b]$. 
Thus if $\{w_i : 1\le i\le \tau\}$ is the set of points sampled by $\pi$,
$\tau_{A}=\sum_{i=1}^\tau \bbone\{w_i\in A\}$. By Lemma \ref{nicepolicieslemma} we
restrict our attention to $\pi\in\bar{\Pi}_{[a,b]}$. In particular if $(x,Y(x))\in K$, $\pi$ will
not choose to sample at $x$ again given initial state $K$. 
Thus 
$\tau_{[a,b]} = \tau_{[a,x]} + \tau_{[x,b]}$ conditioned on any initial history
containing $(x,Y(x))$.  
This is because the only point in $[a,x] \cap [x,b]$ is $\{x\}$, and the lone sample
of $\{x\}$ is in the initial history and it is not counted in $\tau_{[a,x]}$ or $\tau_{[x,b]}$.
From the Markov property, it is also clear that $R_{[a,b]}(H_\tau) = R_{[a,x]}(H_\tau)
+ R_{[x, b]}(H_\tau)$ almost surely conditioned on any initial history containing $x$.

Now, fix some $1< i< t+2$, so that $x_i$ is not $a$ or $b$.
Note that
\begin{eqnarray}
\label{firstline}
V_{[a,b]}(H) & = & 
\sup_{\pi\in\bar{\Pi}_{[a,b]}}\bbE^{\pi}\left[R_{[a,b]}(H_\tau)-c\tau_{[a,b]} \mid H\right]\\
\label{secondline}
 & = & 
 \sup_{\pi\in\bar{\Pi}_{[a,b]}}\left(
 \bbE^{\pi}\left[R_{[a,x_{i}]}(H_\tau)-c\tau_{[a,x_{i}]} \mid H\right]+
 \bbE^{\pi}\left[R_{[x_{i},b]}(H_\tau)-c\tau_{[x_{i},b]} \mid H\right]
 \right)\\
 \label{thirdline}
 & \le & 
 \sup_{\pi\in\bar{\Pi}_{[a,b]}} \bbE^{\pi}\left[R_{[a,x_i]}(H_\tau)-c\tau_{[a,x_i]}\mid H\right]+
 \sup_{\sigma\in\bar{\Pi}_{[a,b]}}\bbE^{\sigma}\left[R_{[x_i,b]}(H_\tau)-c\tau_{[x_i,b]}\mid H\right]\\
 \label{fourthline}
 & = & 
 \sup_{\pi\in\bar{\Pi}_{[a,x_i]}} \bbE^{\pi}\left[R_{[a,x_i]}(H_\tau)-c\tau_{[a,x_i]}\mid H\right]+
 \sup_{\sigma\in\bar{\Pi}_{[x_i,b]}}\bbE^{\sigma}\left[R_{[x_i,b]}(H_\tau)-c\tau_{[x_i,b]}\mid H\right]\\
 & = & 
 \label{fifthline}
 V_{[a,x_{i}]}(H)+V_{[x_{i},b]}(H).
\end{eqnarray}
The equality between (\ref{firstline}) and (\ref{secondline}) holds because $x_i\in H$.
The equality between (\ref{thirdline}) and (\ref{fourthline}) holds because
$\bar{\Pi}_{[a,x_i]}\subseteq\bar{\Pi}_{[a,b]}$ and Lemma \ref{nicepolicieslemma}
shows the supremum is achieved in $\bar{\Pi}_{[a,x_i]}$ and similarly for $\bar{\Pi}_{[x_i,b]}$.
The equality between (\ref{fourthline}) and (\ref{fifthline}) holds 
because $\tau_{[a,x_i]}^\pi=\tau^\pi$
for any $\pi\in\bar{\Pi}_{[a,x_i]}$.

We now show that $V_{[a,b]}(H)\ge V_{[a,x_i]}(H) + V_{[x_i, b]}(H)$.
Let $\pi\in\bar{\Pi}_{[a,x_i]}$ and $\sigma\in\bar{\Pi}_{[x_i,b]}$. 
Define the policy $\gamma$ by
\begin{equation}
\gamma(H)=\left\{ 
  \begin{array}{lr}
    \pi(H\cap [a, x_i]), & \mbox{if }\pi(H\cap [a, x_i])\ne\triangle,\\
    \sigma(H\cap [x_i, b]), & \mbox{otherwise.}
  \end{array}
\right.
\end{equation}

That is, $\gamma$ is the policy that executes $\pi$ with input from $[a,x_i]$ until $\pi$ 
chooses to stop sampling,
and then executes $\sigma$ with input from $[x_i, b]$ until $\sigma$ 
chooses to stop sampling.
Since $\pi\in\bar{\Pi}_{[a,x_i]}$ and $\sigma\in\bar{\Pi}_{[x_i,b]}$ we know
$\tau^{\pi}$ and $\tau^{\sigma}$ are almost surely finite, so $\gamma$ will fully
execute both $\pi$ and $\sigma$. Observe the expectated performance under
$\gamma$ is
\begin{eqnarray}
\bbE\left[R_{[a,b]}(H_\tau^\gamma)-c\tau^\gamma_{[a,b]}\mid H\right] & = & 
\bbE\left[R_{[a,x_{i}]}(H_\tau^\gamma)-c\tau^\gamma_{[a,x_{i}]}\mid H\right]+
\bbE\left[R_{[x_{i},b]}(H_\tau^\gamma)-c\tau^\gamma_{[x_{i},b]}\mid H\right]\\
 & = & 
\bbE\left[R_{[a,x_{i}]}(H_\tau^\pi)-c\tau^\pi_{[a,x_{i}]}\mid H\right]+
\bbE\left[R_{[x_{i},b]}(H_\tau^\sigma)-c\tau^\sigma_{[x_{i},b]}\mid H\right].
\end{eqnarray}
where the decomposition $\tau^\gamma_{[a,b]} = \tau^\gamma_{[a,x_i]} + \tau^\gamma_{[x_i, b]}$ holds 
under $\gamma$ because $(x_i, Y(x_i))\in H$ and so $\gamma$ never samples $x_i$. 
Thus $V_{[a,b]}(H)\ge V_{[a,x_{i}]}(H)+V_{[x_{i},b]}(H)$.
As we have already established the reverse inequality, it follows that
$V_{[a,b]}(H)=V_{[a,x_{i}]}(H)+V_{[x_{i},b]}(H)$.

Now, we partition $H$ about $x_i$: Let $H_{\le i} = \{(x,w)\in H : x\le x_i\}$
and $H_{\ge i} = \{ (x,w) \in H : x \ge x_i \}$.
By Lemma \ref{valueoutsideinterval}, 
$V_{[a,x_{i}]}(H)=V_{[a,x_{i}]}(H_{\le i})$ and 
$V_{[x_{i},b]}(H)=V_{[x_{i},b]}(H_{\ge i})$.
Hence $V_{[a,b]}(H) = V_{[a,x_i]}(H_{\le i}) + V_{[x_i, b]}(H_{\ge i}).$
However, since $x_i$ was chosen to not be an endpoint, $|H_{\le i}| < t+2$
and $|H_{\ge i}| < t+2$. Thus by the induction hypothesis,
\begin{eqnarray}
V_{[a,b]}(H) & = & 
 \sum_{j=1}^{i-1}V_{[x_{j},x_{j+1}]}(H_j)+
 \sum_{j=i}^{t}V_{[x_{j},x_{j+1}]}(H_j)\\
 & = & \sum_{j=1}^{t}V_{[x_{j},x_{j+1}]}(H_j),
\end{eqnarray}
and the induction holds.
\end{proof}

\begin{algorithm}
\caption{Algorithm for computing the value function.
Note the computation on line \ref{theline} is possible because $V[w_1, w_2, x'']$ 
will already be stored for all $w_1,w_2\in\{y_1,\dots,y_n\}$ and $x''\in\{x_1,\dots,x'\}$.
We assume that $Y$ satisfies translation invariance, and thus Proposition \ref{translationthm} applies.}
\begin{algorithmic}[1]
\label{valuealg}
\REQUIRE Interval length $\ell$, $Y$-range discretization $y_1,\dots,y_n$,
$[0,\ell]$-domain discretization $x_1,\dots,x_m$.\\
\ENSURE $x_1=0$ and $x_m=\ell$.
\FOR {$y_L=y_1,\dots,y_n$}
\FOR {$y_R=y_1,\dots,y_n$}
\FOR {$x=x_1,\dots,x_m$}

\IF {$x=0$}
\STATE $V[y_L,y_R,x]\gets 0$
\ELSE
\STATE Let $H=\{(0,y_L),(x,y_R)\}$
\STATE $V[y_L, y_R, x] \gets 
		\max\ \{\bbE[V_{[0,x]}(H\cup \{(x',Y(x')\}) \mid H] - c:x'=x_1,\dots,x\}\cup 
				\{R_{[0,x]}(H)\}$
\label{theline}

\ENDIF

\ENDFOR
\ENDFOR
\ENDFOR

\STATE \textbf{return} $V$

\end{algorithmic}
\end{algorithm}

Theorem \ref{decompositiontheorem} and Proposition \ref{translationthm} give rise to
an efficient algorithm for computing the value function, summarized in Algorithm \ref{valuealg}.
Algorithm \ref{valuealg} takes a discretization $\{x_1,\dots,x_m\}$ of the
domain of $Y$ and discretization $\{y_1,\dots,y_n\}$ of the range of $Y$ as parameters.
It returns the $3$-dimensional array $V[y_L,y_R,x]$ over $y_L,y_R\in\{y_1,\dots,y_n\}$
and $x\in\{x_1,\dots,x_m\}$.  Each element $V[y_L,y_R,x]$ of $V$ is an approximation to $V_{[0,x]}(H)$
where $H=\{(0,y_L), (x,y_R)\}$. Recall Theorem \ref{decompositiontheorem} establishes
that $V_{[a,b]}(H)$ is completely determined by its values on $\{H\in\calH : |H|=2\}$
and by assuming translation invariance Proposition \ref{translationthm} establishes
that if $|H|=2$ then we can assume the leftmost observation in $H$ is at $0$. Thus
the information in the array $V$ can be used to approximate $V_{[a,b]}(H)$ for any $H$.
(Note that finer discretizations lead to more accurate approximations). If $Y$ does not
satisfy translation invariance only a small modification to Algorithm \ref{valuealg} is needed:
$V$ would have to be a $4$-dimensional array, adding one more dimension for the leftmost observation.

The crux of the computation appears on line \ref{theline}. Theorem \ref{decompositiontheorem}
establishes that
\[
V_{[0,x]}\left(\{(0,y_L),(x,y_R),(x',Y(x'))\}\right) = 
V_{[0,x']}\left(\{(0,y_L),(x',Y(x'))\}\right) +
V_{[x',x]}\left(\{(x',Y(x')),(x,y_R)\}\right).
\]
The expectation of the quantities on right side of the above formula will have already
been stored by the algorithm for each point in the $Y$-range discretization, and so
the expectation of the left hand side can be estimated by summing over the $Y$-range 
discretization.

\section{Upper Bound on the Budget-Constrained Problem}
\label{upperboundsec}

So far in this paper we have considered the cost-per-sample scenario,
where the policy may choose how many samples to make without any
additional constraints.  In this section, we show how the cost-per-sample 
problem relates to the budget-constrained problem, in which
the number of samples the policy can take is constrained.  

We first introduce the notion of a randomized policy. Let
$\Pi_R = \{\pi:[0,1]\times\calH \to \bbR_+ \cup \{\Delta\}\}$,
that is, the set of policies which take an additional argument
inside $[0,1]$.  For such policies, we adopt the convention
that histories are still updated according to (\ref{updateequation}),
with the modification that $x_{t+1}=\pi(U, H_t)$ where
$U\sim\mathrm{Uniform}([0,1])$ is drawn once at time $0$ and held fixed over time.
We call these randomized policies because they may take different
actions depending on the random variable $U$.  We will often
write $\pi(H)$ instead of $\pi(U,H)$ when it is clear that
$\pi$ is randomized. Note that taking the supremum in equation
($\ref{valuefunctiondef}$) over the larger set $\Pi_R$ instead of $\Pi$ does not
affect the optimal value, because the deterministic $\epsilon$-optimal policies 
we construct based on Theorem \ref{decompositiontheorem}
remain $\epsilon$-optimal.

Throughout this section we hold an interval $[a,b]$ fixed.
At any state $H\in\calH$ and $T>0$, we define the following sets of constrained policies:
\begin{align}
\Pi_1(H,T) &= \left\{\pi\in\Pi_R : \bbE^\pi \left[\tau\mid H\right] =  T\right\}
&& \text{and} &
\Pi_2(H,T) &= \left\{\pi\in\Pi_R : \bbP^\pi\left(\tau = T\mid H\right)=1\right\}.
\end{align}
The policies in $\Pi_1$ are referred to as expected-budget-constrained policies
and the policies in $\Pi_2$ are referred to as the set of budget-constrained policies.
The corresponding value functions are defined as:
\begin{align}
V_1(H,T) &= \sup_{\pi\in\Pi_1(H,T)} \bbE^\pi \left[ R_{[a,b]}(H_\tau) \mid H\right]
&& \text{and} &
V_2(H,T) &= \sup_{\pi\in\Pi_2(H,T)} \bbE^\pi \left[ R_{[a,b]}(H_\tau) \mid H\right].
\end{align}
Budget-constrained-policies are common in practice - it is
sometimes easier to allocate a predetermined number of samples than to determine
a suitable cost, as the cost-per-sample case requires.
Note the above are defined without a cost. 
This is because $\bbE^\pi\left[\tau\mid H\right]=T$ for any $\pi\in\Pi_1(H,T) \cup \Pi_2(H,T)$,
so any cost term would be constant and not affect the optimal solution.

For the rest of this section we will write the cost-per-sample value function,
as defined in equation (\ref{valuefunctiondef}),
as a function of both the state and the cost.  That is, let $V(H,\lambda)$
indicate $V_{[a,b]}(H)$ with a cost of $\lambda$. Now observe that for any $H\in\calH$,
$\lambda>0$, $T>0$,
\begin{eqnarray}
\label{lineone}
V_1(H,T) & = & \sup_{\pi\in\Pi_1(H,T)} \bbE^\pi 
\left[R_{[a,b]}(H_\tau) -\lambda(\tau - T) \mid H \right] \\
\label{linetwo}
& \le & \sup_{\pi\in\Pi_R} \bbE^\pi \left[R_{[a,b]}(H_\tau) - \lambda \tau \mid H \right] + \lambda T \\
\label{linethree}
& = & \sup_{\pi\in\Pi} \bbE^\pi \left[R_{[a,b]}(H_\tau) - \lambda \tau \mid H \right] + \lambda T \\
& = &  V(H,\lambda) + \lambda T.
\end{eqnarray}
where the inequality between (\ref{lineone}) and (\ref{linetwo}) holds because 
$\Pi_1(H,T)\subseteq\Pi_R$, and the inequality between (\ref{linetwo}) and (\ref{linethree}) holds
because the supremum is attained by a non-randomized policy.
Thus it follows that for any $H\in\calH$,
\begin{equation}
\label{upperboundeq}
V_2(H,T) \le V_1(H,T) \le \inf_{\lambda} V(H,\lambda) + \lambda T
\end{equation}
where the first inequality holds because $\Pi_2(H,T)\subseteq\Pi_1(H,T)$.
Theorem \ref{tightthm} will establish that the second inequality above is tight,
under appropriate assumtions on $T$.

We now introduce some notation.
For any $H\in\calH$ and $\pi\in\Pi_R$ let 
\begin{align}
r(\pi,H)&=\bbE^\pi[R_{[a,b]}(H_\tau)\mid H]
&& \text{and} &
t(\pi,H)&=\bbE^\pi[\tau\mid H].
\end{align}
For any $\lambda > 0$ and $\epsilon>0$, let
\[
\Pi^*_\epsilon(\lambda) = \left\{\pi\in\Pi_R :\forall H\in\calH, 
\per(\pi,\lambda, H) \ge V_{[a,b]}(H,\lambda) - \epsilon  \right\}
\]
denote the set of randomized policies whose performance at every state in $\calH$ with a cost 
of $\lambda$ is at least $\epsilon$-optimal.
Let
\begin{equation}
\label{slambdadef}
s(\lambda, H)=\limsup_{\epsilon\to0^+}\{t(\pi,H) : \pi\in\Pi^*_\epsilon(\lambda)\}
\end{equation}
denote the limit as $\epsilon$ decreases to $0$ of the maximum expected number of 
samples an $\epsilon$-optimal policy takes at any state $H$ given a
cost of $\lambda$. The function $s(\lambda,H)$ will be instrumental 
in characterizing the policies in $\Pi_1(H,T)$ and thus the expected-budget-constrained
value function $V_1(H,T)$. Finally, we define
\begin{equation}
\bar{T}(H)=\lim_{\lambda\to 0^+} s(\lambda,H).
\end{equation}
Intuitively, $\bar{T}(H)$ denotes the maximum number of samples a sensible policy
takes as the cost of taking samples decreases to $0$. Note that $\bar{T}(H)$ inherently
depends on the Markov process $Y$. It is natural to expect that as the cost decreases to $0$
optimal policies will begin to take more and more samples, and hence $\bar{T}(H)=\infty$.
However, it is possible to construct Markov processes $Y$ that are completely characterized
by finitely many samples within $[a,b]$. 

The main result of this section is stated in the following theorem.

\begin{thm}
\label{tightthm}
Fix a state $H_0\in\calH$ and let $0\le T < \bar{T}(H_0)$. Then
\[
V_1(H_0,T) = \inf_\lambda V(H_0,\lambda) + \lambda T.
\]
\end{thm}

Theorem \ref{tightthm} is nice because it gives an alternate
characterization of $V_1(H,T)$, but its main importance comes in noting the following:
\begin{eqnarray*}
V(H,\lambda) + \lambda T & = & 
\sup_{\pi\in\Pi}r(\pi,H) - \lambda t(\pi,H) + \lambda T\\
& = & 
\sup_{\pi\in\Pi}r(\pi,H) +\lambda(T -t(\pi,H)).
\end{eqnarray*}
The fact that $V(H,\lambda)$ is a supremum over linear functions of $\lambda$ implies
that it is a convex function of $\lambda$.  Furthermore, for each value of $\lambda$
the quantity $V(H,\lambda)+\lambda T$ can be computed using simulations estimate
the performance of the policy defined in Algorithm \ref{valuealg}.
Thus the actual value of $V_1(H)$ can be computed as the solution to a convex
program in $\lambda$, which can be solved easily by algorithms such as bisection search.

\begin{proof}[\textbf{Proof of Theorem \ref{tightthm}}]
We hold $H_0\in\calH$ fixed throughout the proof, and for simplicity omit $H$ from
the notation defined above (e.g. we refer to $s(\lambda,H_0)$ as $s(\lambda)$).

We first note two properties concerning $s(\lambda)$:
\begin{enumerate}
\item $\lim_{\lambda\to\infty}s(\lambda) = 0$. From inequality (\ref{rewardupperbound})
we know there exists some constant $C>0$ such that $|R_{[a,b]}(H_0)| \le C$. Thus the
policy that takes $0$ samples achieves a performance of $-C$ at worst. Take $\lambda>2C$.
Then the best performance a policy that takes $1$ or more samples can achieve is 
worse than $C-2C=-C$, meaning that the $0$-sample-policy is best for such $\lambda$.
\item 
$s(\lambda)$ is monotonically decreasing. 
Fix $\epsilon>0$, let $\lambda < \lambda'$
and pick any $\pi\in\Pi^*_\epsilon(\lambda)$ and $\pi'_\epsilon\in\Pi^*(\lambda')$. By the 
$\epsilon$-optimality
of $\pi$ with respect to $\lambda$ and $\pi'$ with respect to $\lambda'$, the
following inequalities hold:
\begin{eqnarray}
\label{ineq1}
r(\pi) - \lambda t(\pi) & \ge & r(\pi') - \lambda t(\pi') -\epsilon \\
\label{ineq2}
r(\pi')- \lambda' t(\pi') & \ge & r(\pi) - \lambda' t(\pi) - \epsilon.
\end{eqnarray}
Subtracting (\ref{ineq2}) from (\ref{ineq1}) it follows that 
$(\lambda' - \lambda) t(\pi) +\epsilon \ge (\lambda' - \lambda)t(\pi')-\epsilon$.
Since $\lambda' > \lambda$, it follows that
\begin{equation}
\label{tpiepsilonineq}
t(\pi) \ge t(\pi') - \frac{2\epsilon}{\lambda'-\lambda}.
\end{equation}
Taking $\epsilon\to0^+$, we conclude $s(\lambda)\ge s(\lambda')$.

\end{enumerate}

Since $s(\lambda)$ converges to $\bar{T}(H_0)$ as $\lambda\to 0^+$,
converges to $\infty$ as $\lambda\to\infty$, is monotonically decreasing,
and $0\le T < \bar{T}(H_0)$, it follows that there exists
some $\lambda^*$ such that either $s(\lambda^*)=T$ or there is a jump discontinuity at $\lambda^*$
around $T$, that is $\lim_{\lambda\to\lambda^{*+}}s(\lambda) \le T$ and 
$\lim_{\lambda\to\lambda^{*-}} s(\lambda) \ge T$.
Pick any sequence $(\overlambda_n, \underlambda_n)$ such that
$(\overlambda_n)$ is decreasing in $n$, $(\underlambda_n)$ is
increasing in $n$, and 
$\lim_{n\to\infty}\overlambda_n = \lambda^* = \lim_{n\to\infty}\underlambda_n$.
For each $n\in\bbN$ and $\epsilon>0$, pick any $\overpi_n^\epsilon\in\Pi^*_\epsilon(\overlambda_n)$
and $\underpi_n^\epsilon\in\Pi^*_\epsilon(\underlambda_n)$. 
By definition of $s(\lambda)$ and the fact that $\overlambda_n > \lambda^*$, we know 
\begin{equation}
t(\overpi_n^\epsilon) \le T. 
\end{equation}
Similarly, from equation (\ref{tpiepsilonineq}) and the fact that 
$\underlambda_n < \lambda^*$, we know 
\begin{equation}
t(\underpi_n^\epsilon) \ge T - g_n(\epsilon),
\end{equation}
where $g_n(\epsilon)=\frac{2\epsilon}{\underlambda_n - \lambda^*}$.

We now construct a sequence of randomized policies $\{\pi_n^\epsilon\}_{n\in\bbN}$ 
based on $\{\underpi_n^\epsilon\}_{n\in\bbN}$ and $\{\overpi_n^\epsilon\}_{n\in\bbN}$.
First, we define probabilities $p_n^\epsilon$ by 
\begin{equation}
p_n^\epsilon \coloneqq 
\begin{cases}
\frac{T-t(\overpi_n^\epsilon)}{t(\underpi_n^\epsilon) - t(\overpi_n^\epsilon)}, 
&\text{if $t(\underpi_n^\epsilon) \ne t(\overpi_n^\epsilon)$ and $t(\underpi_n^\epsilon)\ge T$,}\\
\frac{1}{2}, 
&\text{if $t(\underpi_n^\epsilon) = t(\overpi_n^\epsilon)$ and $t(\underpi_n^\epsilon)\ge T$,} \\
1, &\text{otherwise.}
\end{cases}
\end{equation}
The policies are defined by
\begin{equation}
\pi_n^\epsilon(H, U) \coloneqq \begin{cases}
\overpi_n^\epsilon(H),&\text{if $U > p_n^\epsilon$}\\
\underpi_n^\epsilon(H), &\text{if $U \le p_n^\epsilon$.}
\end{cases}
\end{equation}
Note that all of these policies satisfy
$T\ge t(\pi_n^\epsilon)\ge T - g_n(\epsilon)$.

We now turn to a technical equality involving the policies $\{\underpi_n^\epsilon\}$
and $\{\overpi_n^\epsilon\}$. Let
\begin{equation}
L_n^\epsilon=p_n^\epsilon \underlambda_n [t(\underpi_n^\epsilon) - T] +
(1-p_n^\epsilon)\overlambda_n [t(\overpi_n^\epsilon) - T].
\end{equation}
We claim the following equality holds:
\begin{equation}
\label{technicalineq}
\lim_{n\to\infty} \liminf_{\epsilon\to0^+}
L_n^\epsilon = 0.
\end{equation}
To see this, first note that $L_n^\epsilon$ can be rewritten as
\begin{eqnarray}
\label{firstlinesec4}
L_n^\epsilon & = &
p_n^\epsilon\overlambda_n [t(\underpi_n^\epsilon)- T] + 
(1-p_n^\epsilon)\overlambda_n [t(\overpi_n^\epsilon) - T]
- p_n^\epsilon (\overlambda_n - \underlambda_n) [t(\underpi_n^\epsilon) - T] \\
\label{secondlinesec4}
& = & \overlambda_n[t(\pi_n^\epsilon) - T] - 
p_n^\epsilon(\overlambda_n - \underlambda_n) [t(\underpi_n^\epsilon) - T].
\end{eqnarray}
The equality (\ref{firstlinesec4}) can be derived from simple algebra, and 
equality (\ref{secondlinesec4}) holds due to the definition of $\pi_n^\epsilon$.
Since $T\ge t(\pi_n^\epsilon) \ge T - g_n(\epsilon)$ and $\lim_{\epsilon\to 0^+}g_n(\epsilon)=0$
for each $n\in\bbN$, it follows that 
\begin{equation}
\liminf_{\epsilon\to 0^+}\overlambda_n[t(\pi_n^\epsilon) - T] = 0.
\end{equation}
Similarly, note that 
\begin{equation}
\liminf_{\epsilon\to 0^+}p_n^\epsilon(\overlambda_n - \underlambda_n)[t(\underpi_n^\epsilon) - T]
= (\overlambda_n - \underlambda_n)
\liminf_{\epsilon\to 0^+} p_n^\epsilon[t(\underpi_n^\epsilon) - T].
\end{equation}
Since the $(p_n^\epsilon)$ and $[t(\underpi_n^\epsilon)-T]$ are bounded in $\epsilon$
the $\liminf$ above is finite. Since $\lim_{n\to\infty}\overlambda_n - \underlambda_n = 0$,
\begin{equation}
\lim_{n\to\infty}\liminf_{\epsilon\to 0^+}L_n^\epsilon = 0 - 
\lim_{n\to\infty}(\overlambda_n-\underlambda_n)
\liminf_{\epsilon\to 0^+} p_n^\epsilon[t(\underpi_n^\epsilon) - T]
=0
\end{equation}
and (\ref{technicalineq}) is established.

Now let $S=\inf_{\lambda} V(H,\lambda) + \lambda T$.
For any $n\in\bbN$ and $\epsilon>0$, because $S\le V(H,\underlambda_n) + \underlambda_n T$ and similarly
for $\overlambda_n$, we have
\begin{eqnarray}
S & \le &
p_n^\epsilon \left[V(H,\underlambda_n) + \underlambda_n T \right] +
(1-p_n^\epsilon)\left[V(H,\overlambda_n) + \overlambda_n T \right] \\
& \le &
p_n^\epsilon \left[r(\underpi_n^\epsilon) - 
\underlambda_n (t(\underpi_n^\epsilon) - T) + \epsilon\right] + 
(1-p_n^\epsilon) \left[r(\overpi_n^\epsilon) - 
\overlambda_n(t(\overpi_n^\epsilon) - T) + \epsilon\right] \\
& = &
r(\pi_n^\epsilon) - p_n^\epsilon \underlambda_n [t(\underpi_n) - T] -
(1-p_n^\epsilon)\overlambda_n [t(\overpi_n) - T] +\epsilon.
\end{eqnarray}
Taking $\liminf_{\epsilon\to 0^+}$ on both sides, it follows that
\begin{eqnarray}
S & \le & \liminf_{\epsilon\to 0^+} \left(
r(\pi_n^\epsilon) - p_n^\epsilon \underlambda_n [t(\underpi_n) - T] -
(1-p_n^\epsilon)\overlambda_n [t(\overpi_n) - T] +\epsilon\right) \\
& \le & 
\liminf_{\epsilon\to 0^+} r(\pi_n^\epsilon) - 
\liminf_{\epsilon\to 0^+} \left(p_n^\epsilon \underlambda_n [t(\underpi_n^\epsilon) - T] +
(1-p_n^\epsilon)\overlambda_n [t(\overpi_n^\epsilon) - T]\right).
\end{eqnarray}
Equation (\ref{technicalineq}) establishes that taking $n\to\infty$ 
sends the second above term to $0$. Hence
\begin{equation}
S\le\lim_{n\to\infty}\liminf_{\epsilon\to 0^+} r(\pi_n^\epsilon).
\end{equation}

We now note that for any policy $\pi\in\Pi_R$ such that $t(\pi)\le T$,
there exists a policy $\pi'\in \Pi_1(T)$ such that $r(\pi)=r(\pi')$.
Define the policy $\pi'$ by letting $\pi'(H)=\pi(H)$ provided $\pi(H)\ne\Delta$. 
Let $h\in H_0$ be any point that has already been samples. 
Once $\pi(H)=\Delta$, $\pi'$ chooses to sample at $h$ for $\lfloor{T-t(\pi)}\rfloor$ iterations.
Finally, $\pi'$ chooses to sample at $h$ one more time with probability
$\lfloor{T - t(\pi)}\rfloor-(T-t(\pi))$. By construction it follows that $t(\pi')=t(\pi)$,
and since sampling at $h$ does not affect the reward we have that $r(\pi')=r(\pi)$.

Since $T\ge t(\pi_n^\epsilon)$ for all $n$ and $\epsilon$, it follows that 
$V_1(T)\ge \sup_{n,\epsilon} t(\pi_n^\epsilon)$. Thus,
\begin{eqnarray*}
\sup_{n,\epsilon} r(\pi_n^\epsilon) & \le & V_1(T) \\
& \le & S \\
& \le & \lim_{n\to\infty}\liminf_{\epsilon\to 0^+} r(\pi_n^\epsilon) \\
& \le & \sup_{n,\epsilon} r(\pi_n^\epsilon).
\end{eqnarray*}
Since the first and final terms above are the same, the inequalities above are forced to be equalities,
and so $S=V_1(T)$.
\end{proof}

\section{Experimental Analysis}
\label{simulationsec}

In this section we run simulations to better understand the behaviour
of the optimal policy. 
We focus on superlevel set detection, and consider two 
choices for the Markov process $Y$.  The first
is a standard Brownian motion.  The second is a compound Poisson
process, that is
\[
Y(t) = \sum_{i=0}^{N(t)} D_i,
\]
where $N(t)$ is a Poisson process with parameter $\mu$ and $D_i$ are
independent standard normal variables. We consider the interval 
$[a,b]=[0,1]$, the threshold $k=0$ and assume we have observed
endpoint observations $Y(0)=Y(1)=0$. For the compound Poisson process
we use a parameter of $\mu=20$.  
The algorithm we use to compute the value function (and thus the optimal policy)
is given in Algorithm \ref{valuealg}. 
Both the standard Brownian motion
and compound Poisson process satisfy translation invariance
(as defined in equation (\ref{translationinvariance})), so we are
only concerned with a $3$-dimensional state space.
To compute the optimal
policy we discretize the domain and range of $Y$. For all of our experiments
we use the indicator reward functions: $f_+(y)=\bbone\{y\ge k\}$ and $f_-(y)=\bbone\{x\le k\}$.

\begin{figure}

\begin{centering}
\includegraphics[scale=0.33]{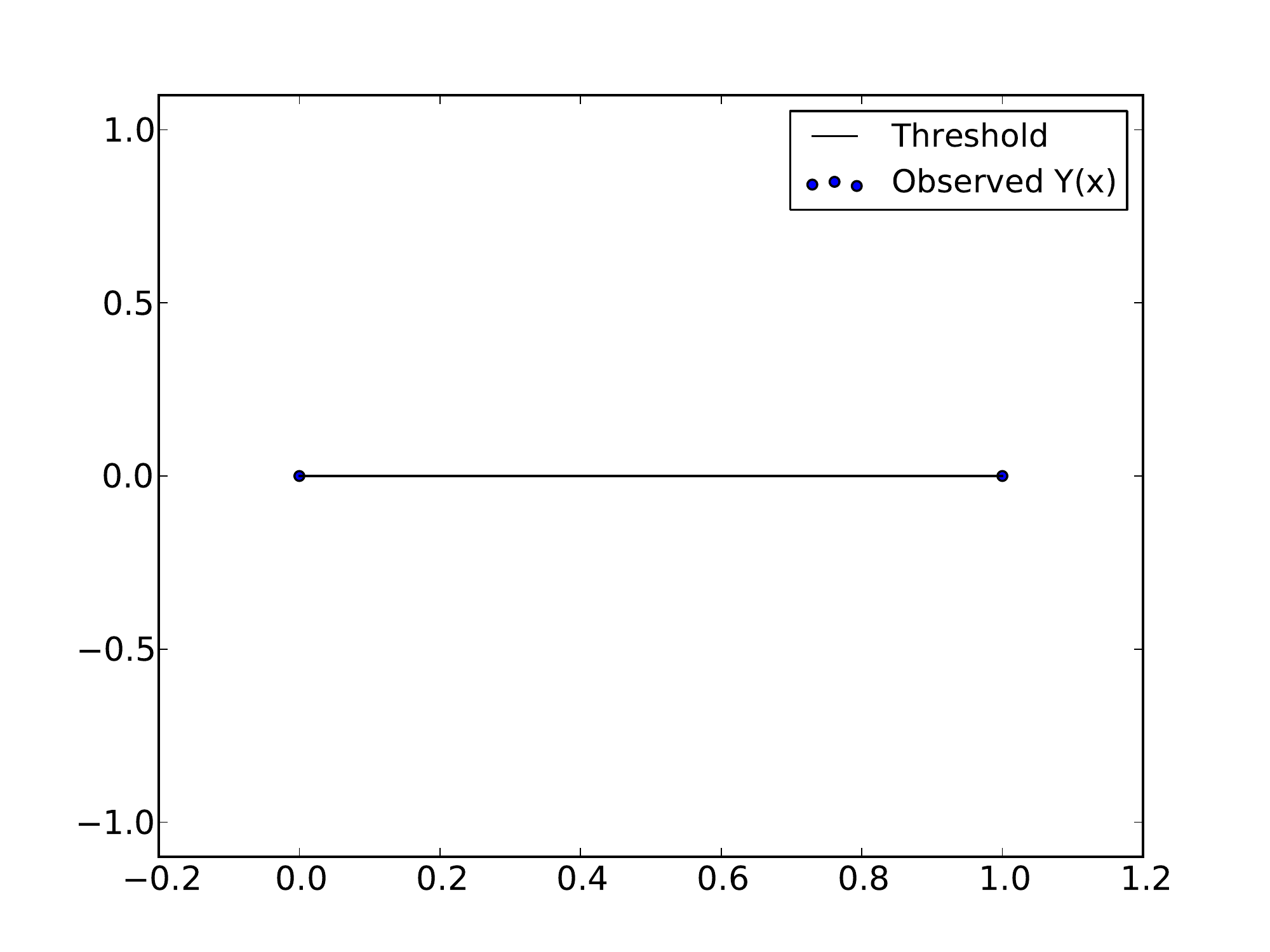}\includegraphics[scale=0.33]{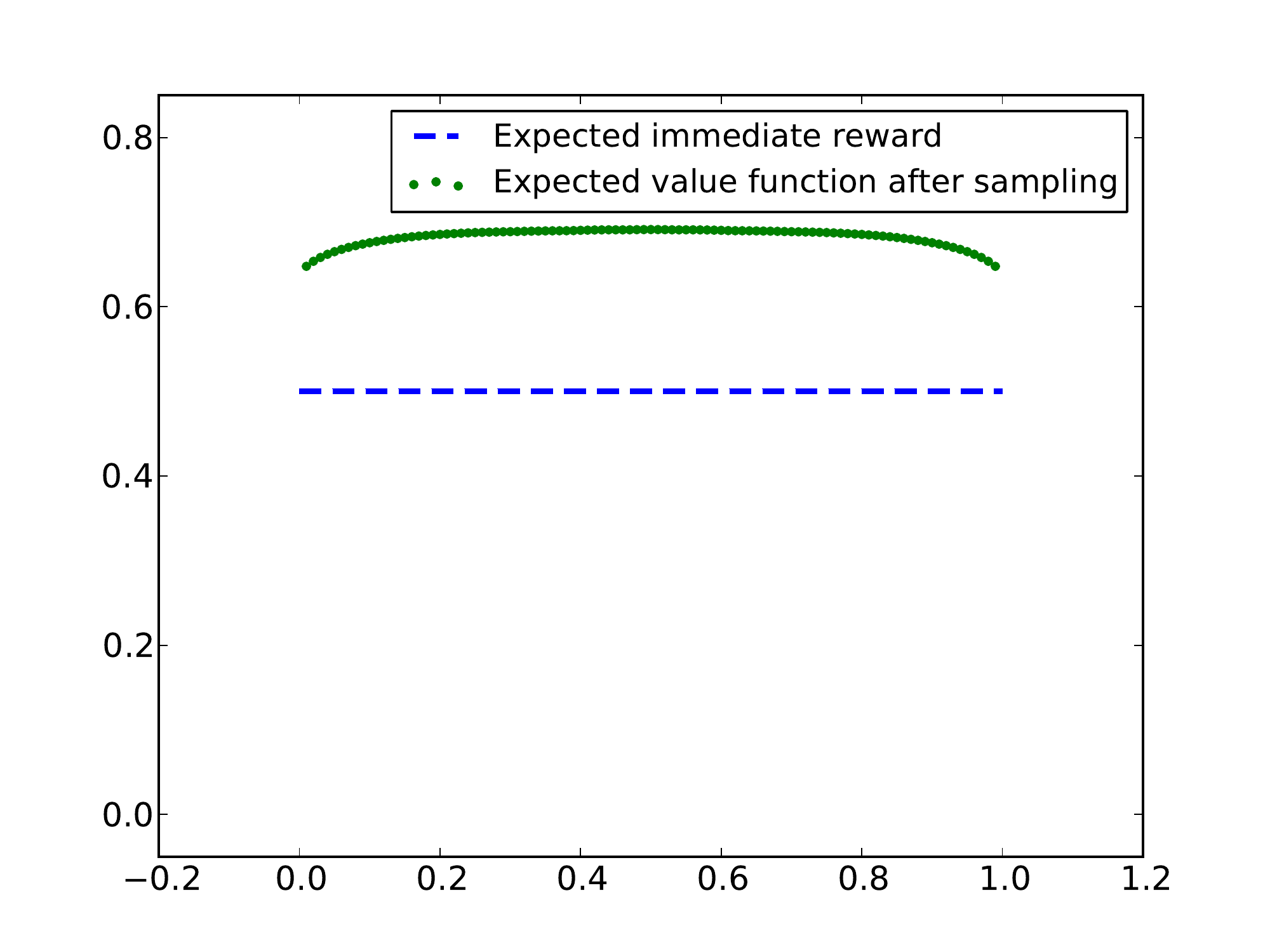}
\par\end{centering}

\begin{centering}
\includegraphics[scale=0.33]{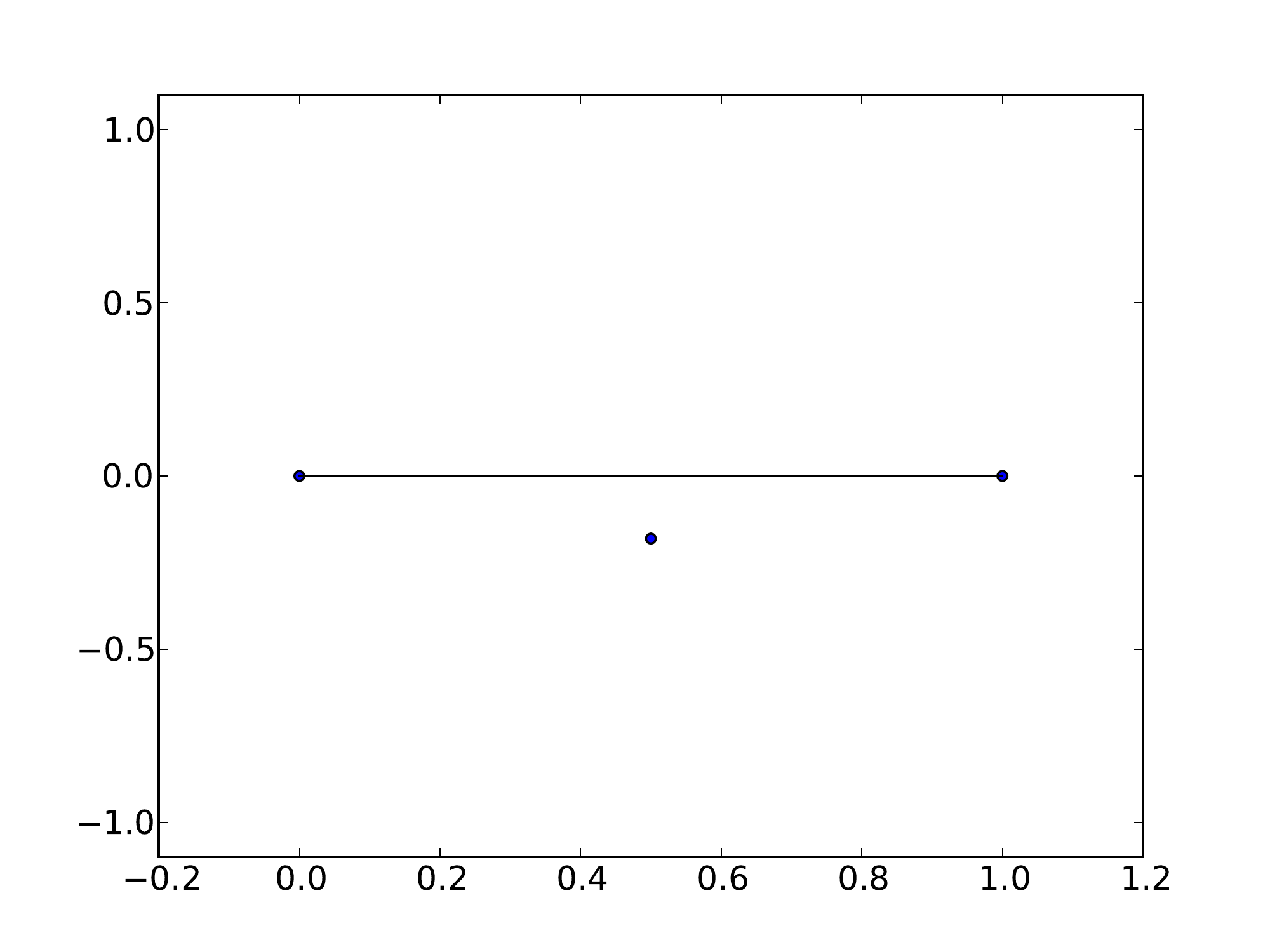}\includegraphics[scale=0.33]{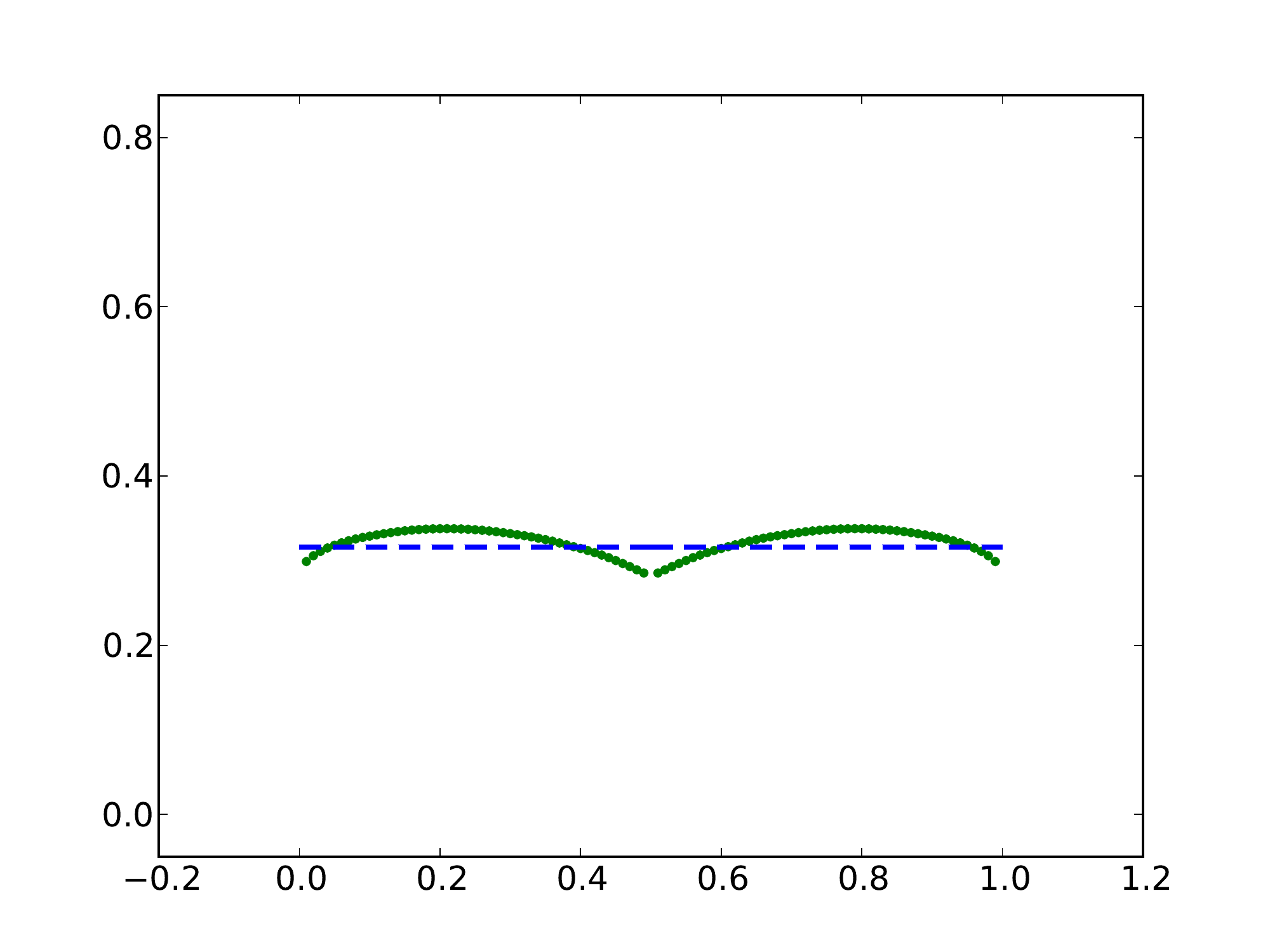}
\par\end{centering}

\begin{centering}
\includegraphics[scale=0.33]{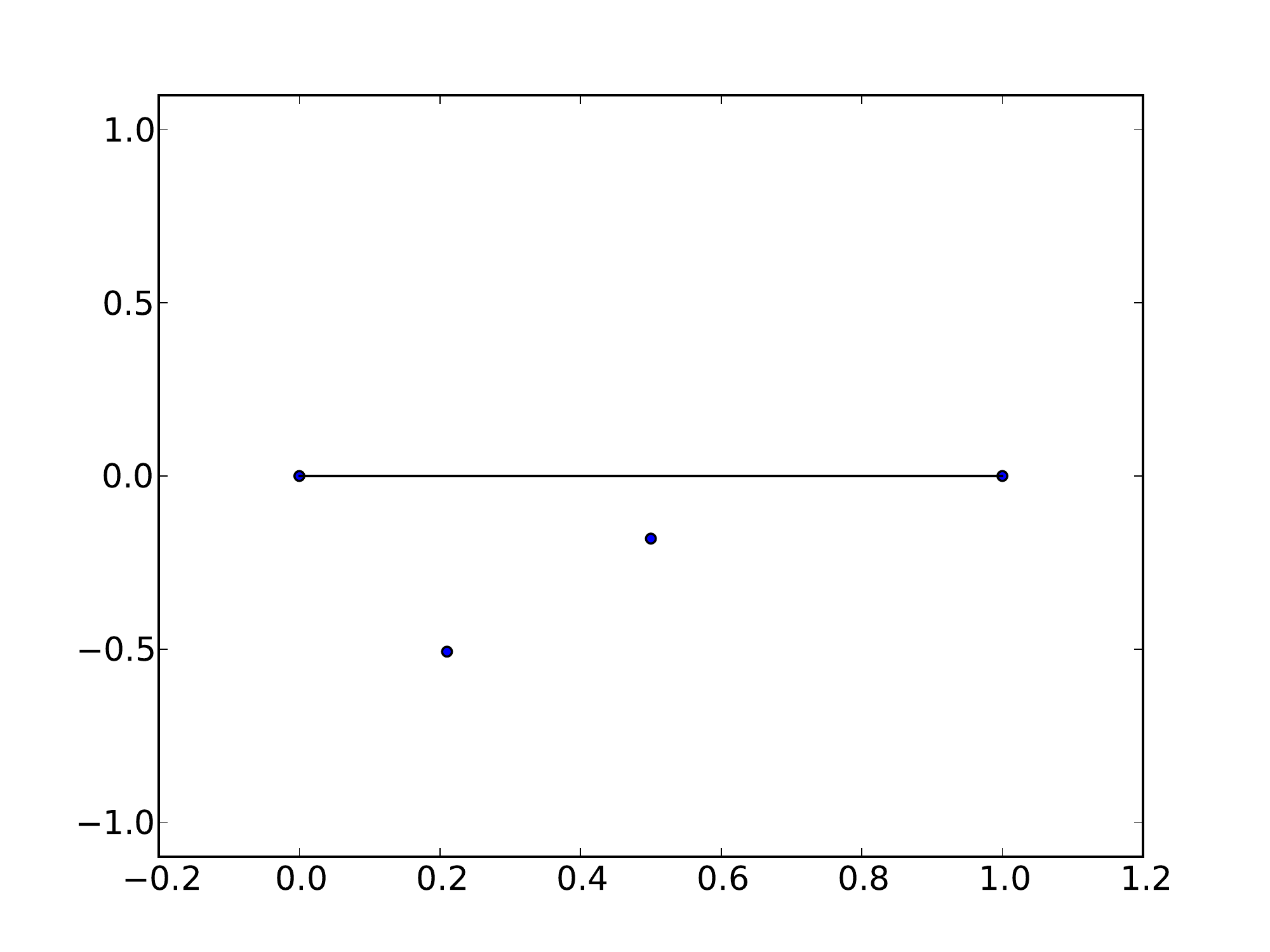}\includegraphics[scale=0.33]{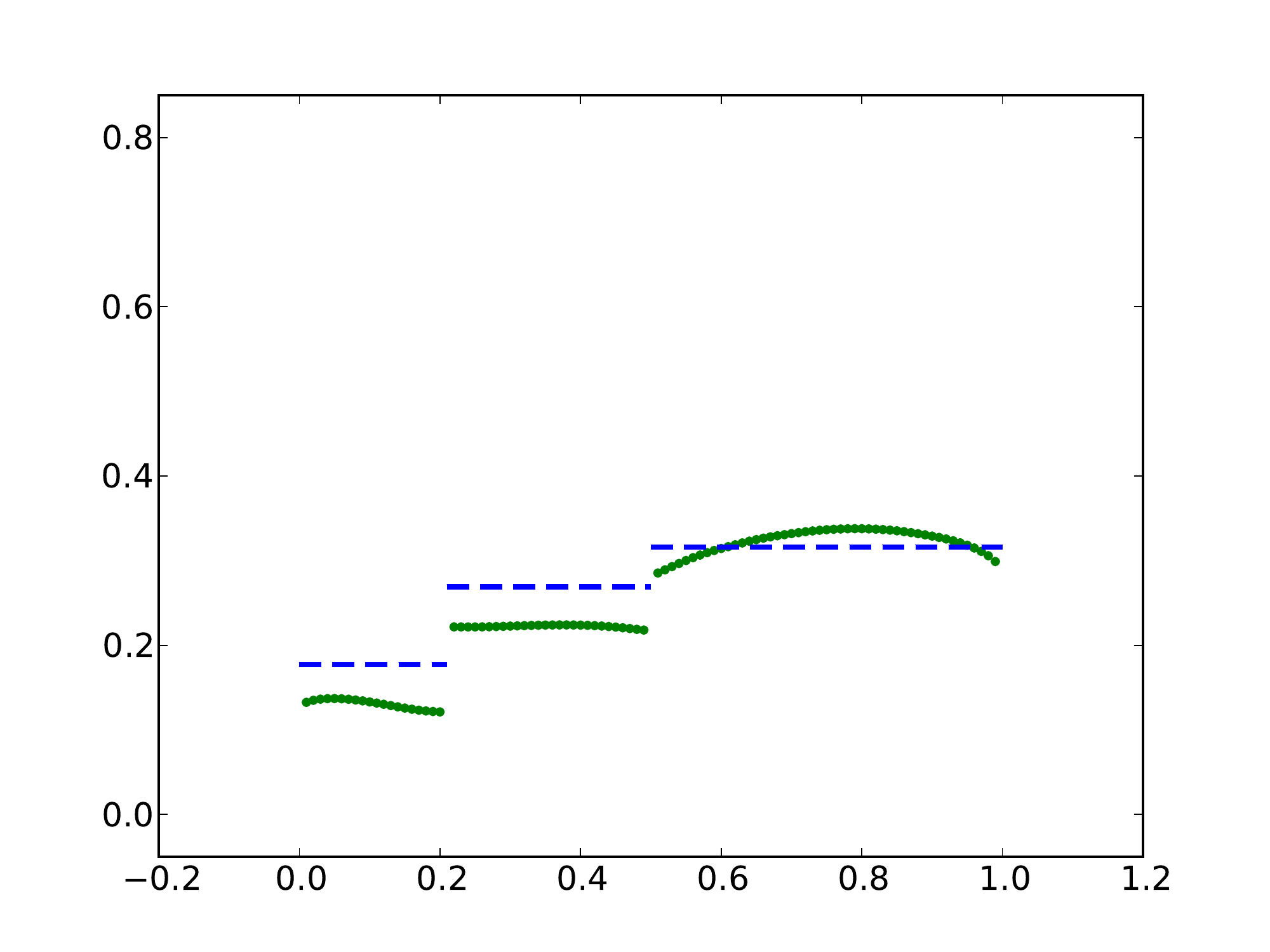}
\par\end{centering}

\begin{centering}
\includegraphics[scale=0.33]{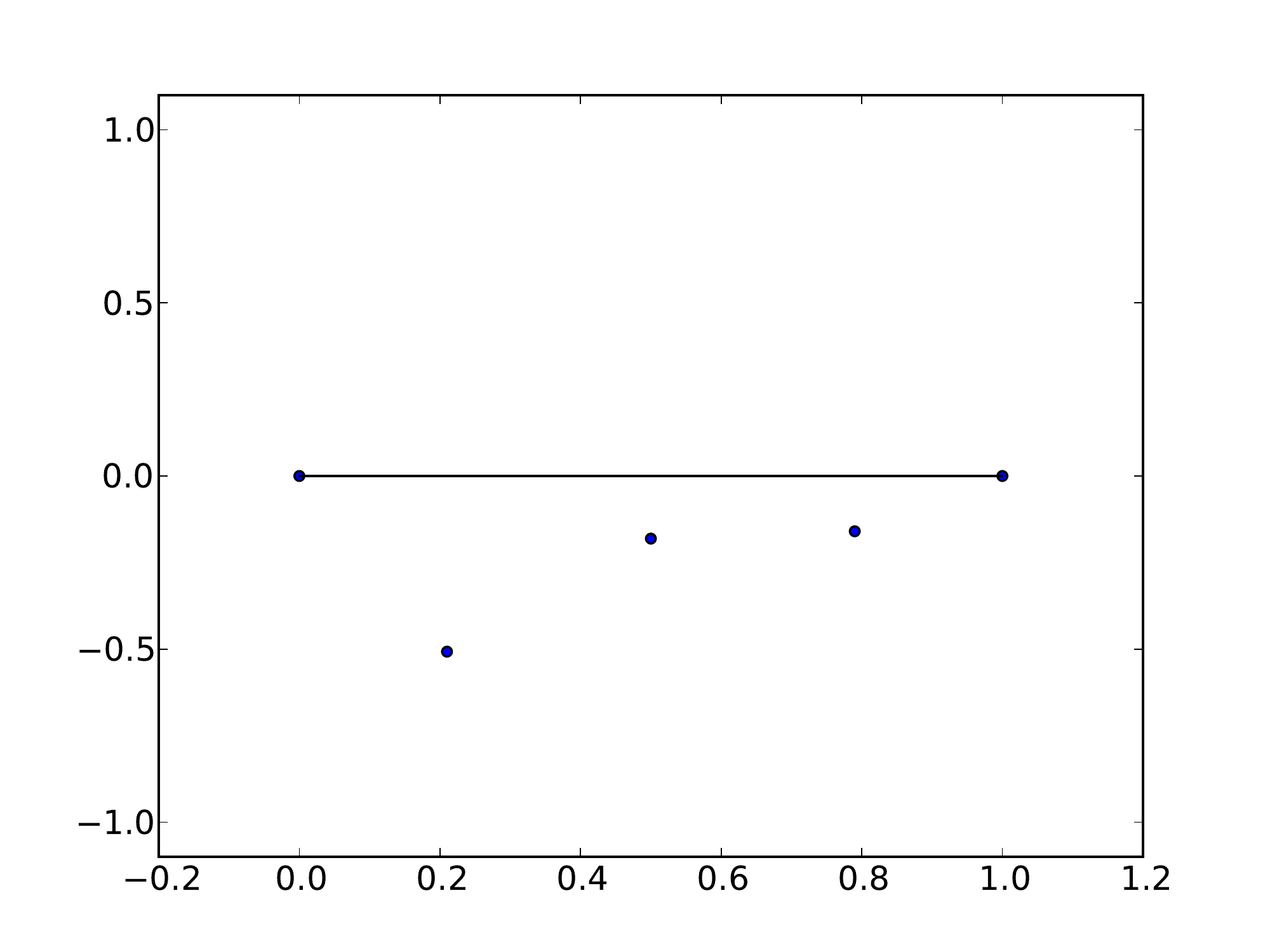}\includegraphics[scale=0.33]{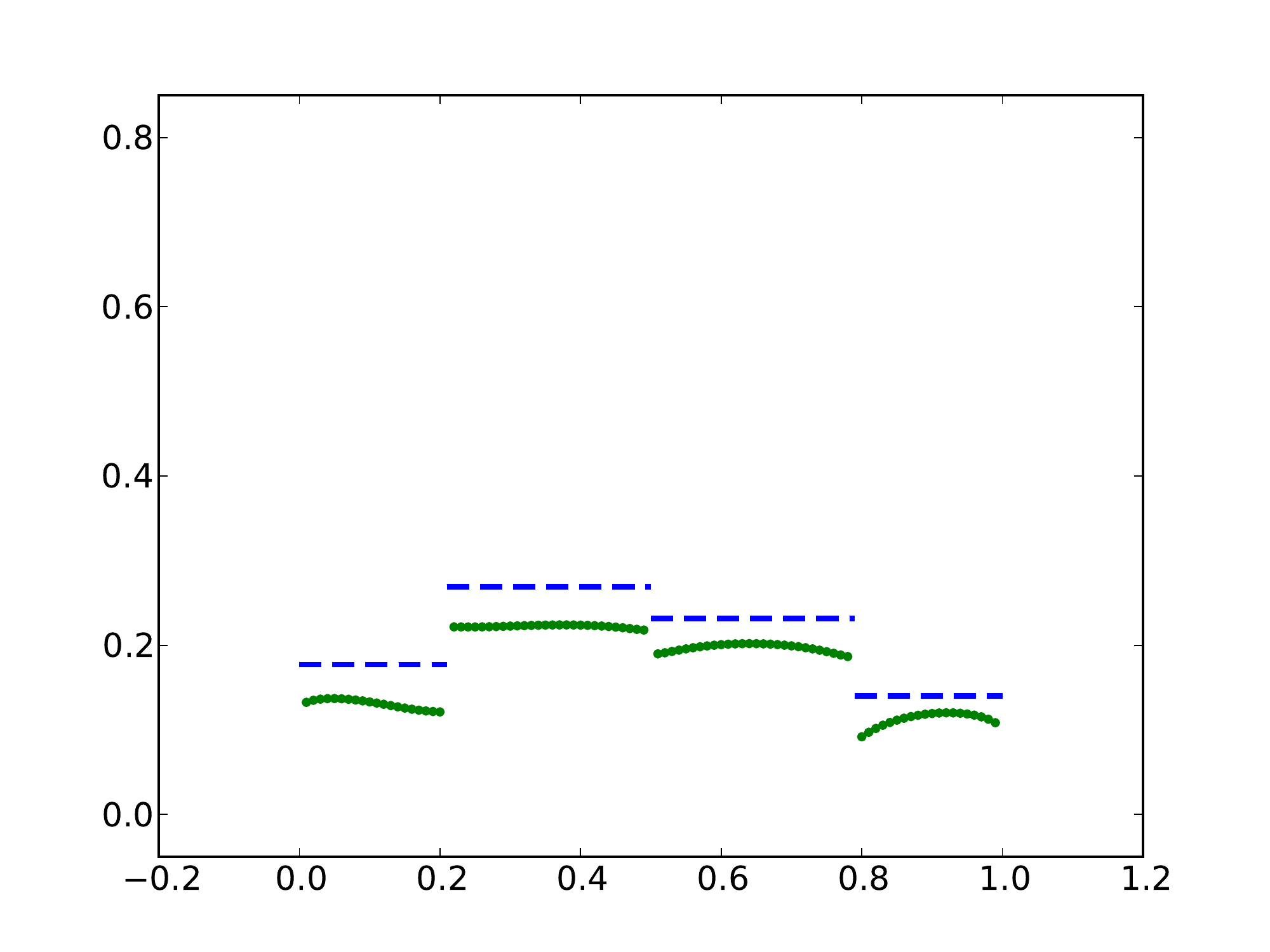}
\par\end{centering}

\caption{Depiction of an optimal policy. 
Here $Y$ is a standard Brownian motion, 
set the threshold $k=0$, x-axis discretization of $100$, cost $c=0.05$. 
\textbf{Left:} Sampled points at each iteration. 
\textbf{Right:} Expected value of the value function (solid line) plotted against expected reward (dashed line). 
The policy samples the point maximizing the difference between these quantities until the maximum difference is negative.
}

\label{samplingexample}
\end{figure}

Figure \ref{samplingexample} depicts the behavior of the optimal policy 
defined in Algorithm \ref{valuealg}
for a Brownian motion over the interval $[0,1]$. 
The optimal policy and expected value function exhibits several intuitive properties. 
First, note that the difference between the expected value of sampling and the
expected reward of not sampling is larger for the intervals where
the endpoints are further away from the threshold. Second, note that the
optimal policy takes its first sample exactly in the middle of the interval $[0,1]$ at $x=0.5$,
the point with the highest variance.

\begin{figure}
\centering
\begin{subfigure}[b]{0.45\textwidth}
\includegraphics[width=\textwidth]{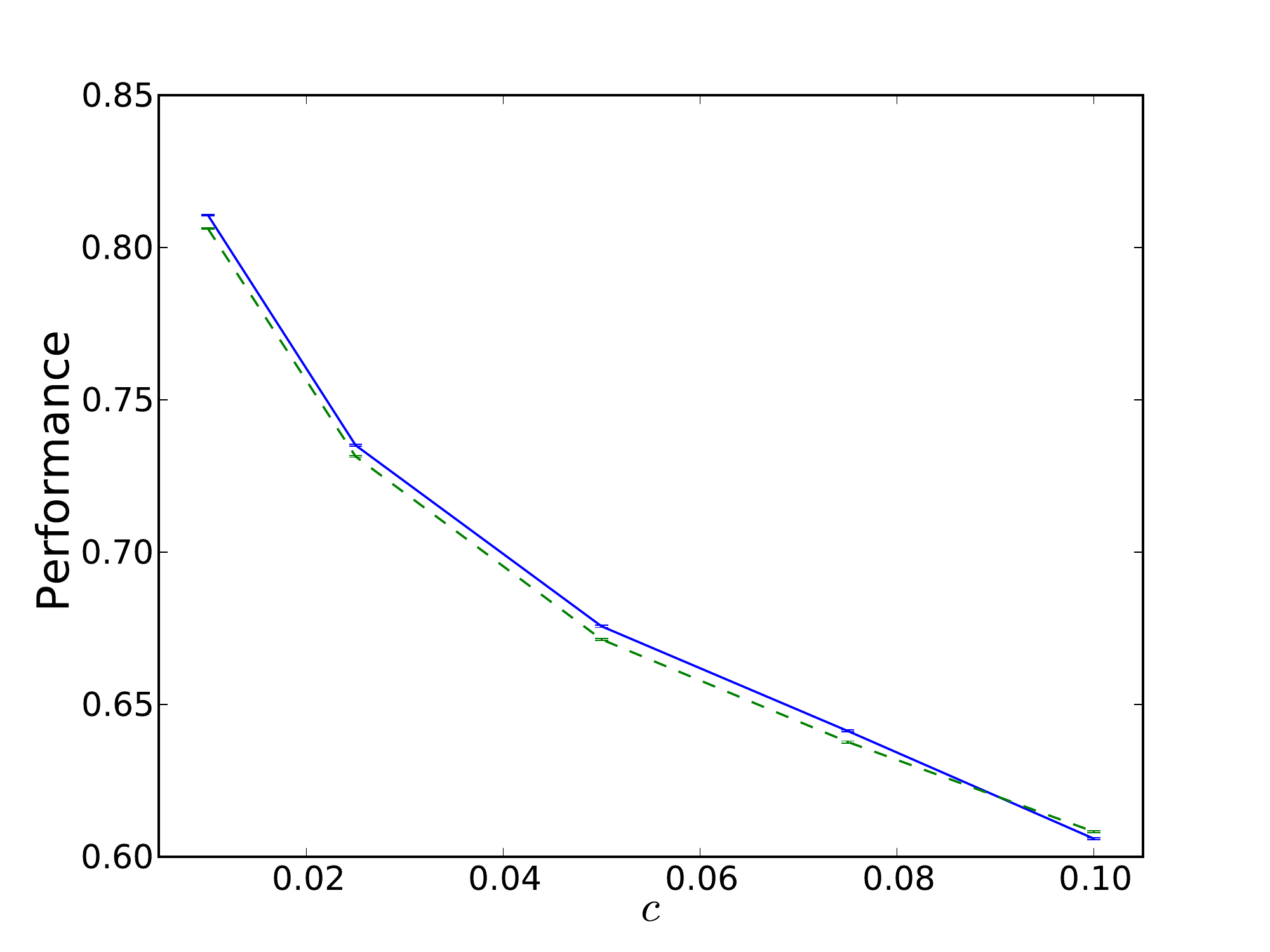}
\caption{$Y$ is a Brownian motion.}
\end{subfigure}
\begin{subfigure}[b]{0.45\textwidth}
\includegraphics[width=\textwidth]{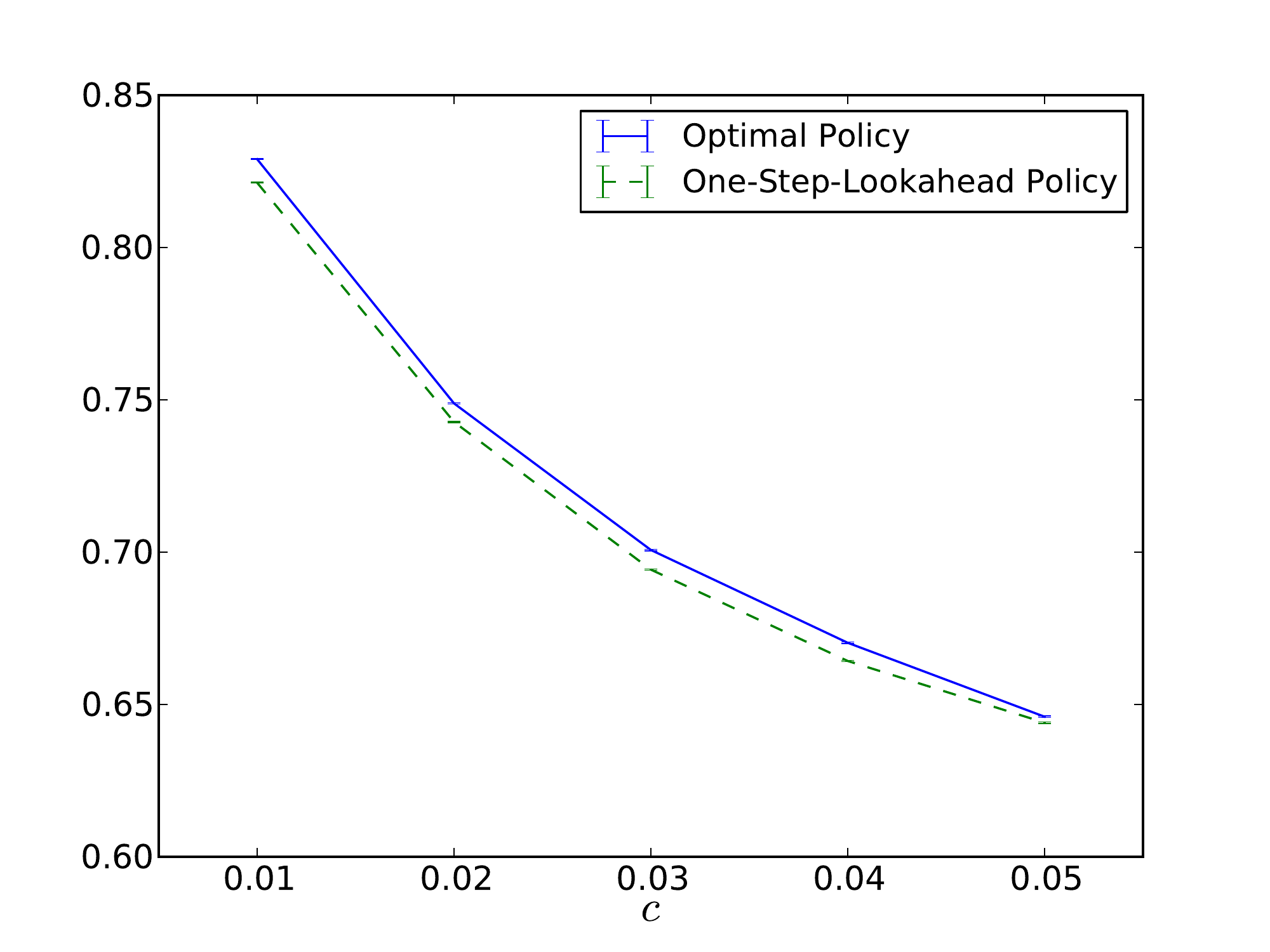}
\caption{$Y$ is a compound Poisson process.}
\end{subfigure}
\caption{Value of optimal policy (in blue) and value of 
one-step lookahead policy (in green) vs. the cost of sampling. }
\label{onestepcomparison}
\end{figure}

One benefit of being able to compute an optimal policy is being able
to characterize suboptimality of the common one-step lookahead heuristic policy
described in the introduction.  
In this problem setup, the one-step lookahead heuristic policy samples
the point maximizing the expected immediate reward, or chooses to stop
sampling when the gain in reward is lower than the cost.
We simulated both this policy and the optimal policy,
varying the cost $c$. We used both the compound Poisson process and Brownian
motion prior on $Y$, again over $[0,1]$. As above, we
assume the  initial history $H_0=\{(0,0), (0,1)\}$, that is, we assume initial 
observations $Y(0)=Y(1)=0$, and use a threshold of $k=0$.
The value of a policy is estimated by running the policy under the above
conditions (where we sample the observations of $Y$ from the corresponding
conditional distribution), and looking at the expected reward once the policy chooses to
stop sampling. We obtain accurate estimates by running each policy $100000$ times.
The results are shown in Figure \ref{onestepcomparison}. While the one-step lookahead
policy is clearly suboptimal, for each value of $c$ we tested, the value of the one-step
lookahead policy is within $98$\% of the optimal.  This bodes well for the performance of more realistic one-step
lookahead algorithms, which are common in practice as discussed in Section \ref{introsection}.

\begin{figure}
\centering
\begin{subfigure}[b]{0.45\textwidth}
\includegraphics[width=\textwidth]{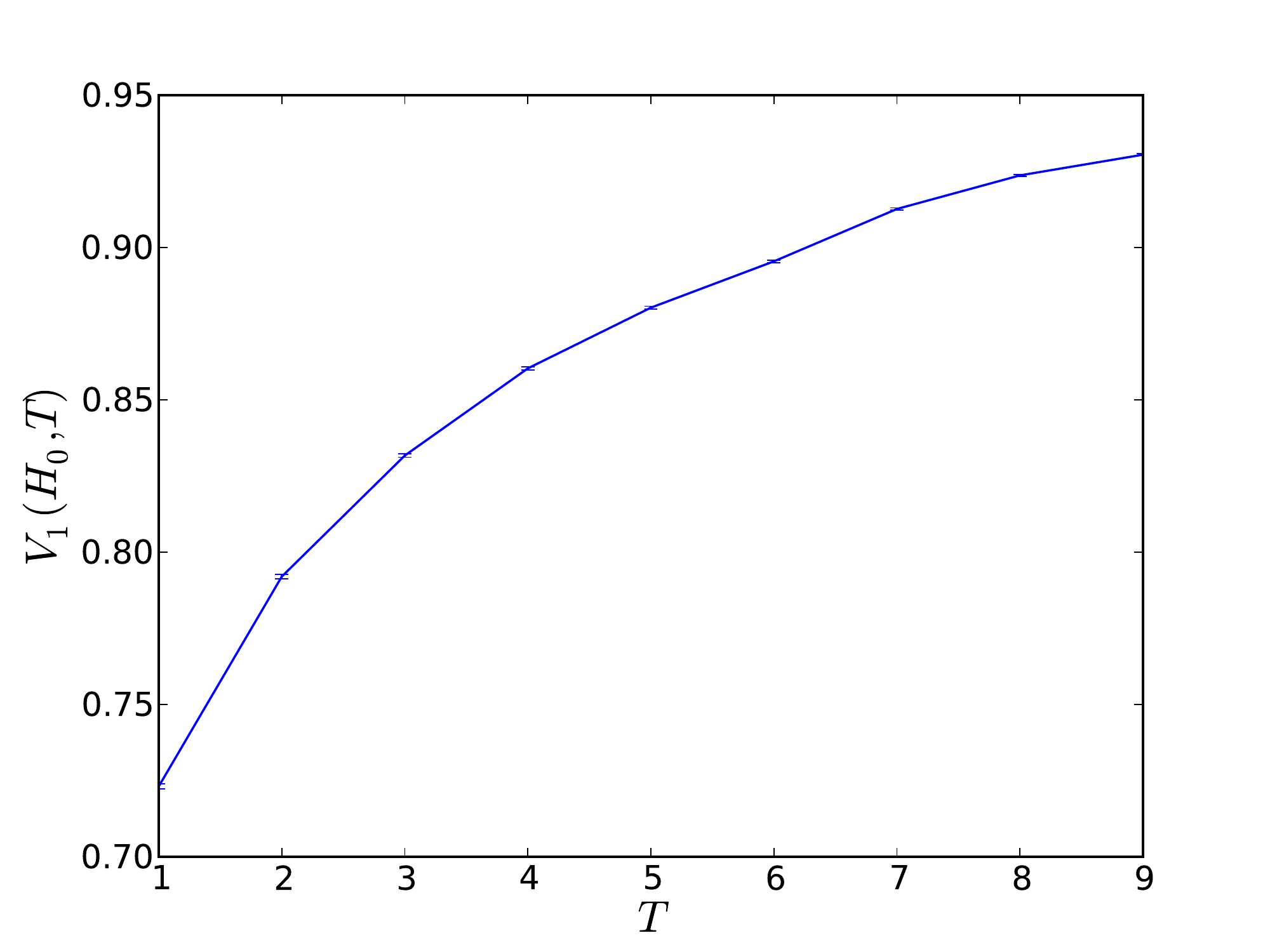}
\caption{Value of expected constrained budget constrained problems.}
\label{expectedconstrainedvalue}
\end{subfigure}
\begin{subfigure}[b]{0.45\textwidth}
\includegraphics[width=\textwidth]{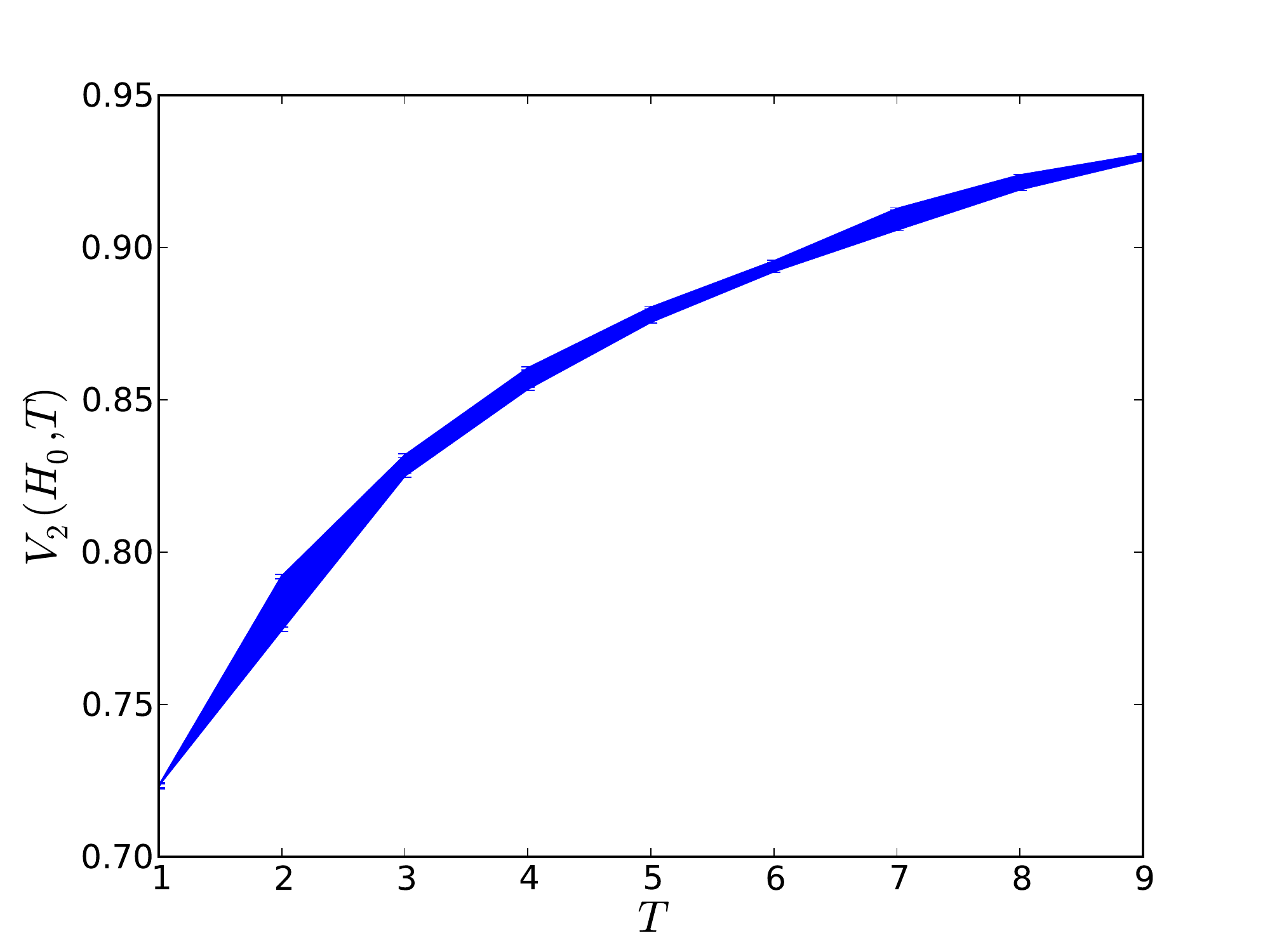}
\caption{Region containing value of exact constrained budget problems.}
\label{tstepvalues}
\end{subfigure}
\caption{Constrained-budget value plots (expected and almost sure contraints on the
left and right respectively) 
when $Y$ is a Brownian motion. On the right, the upper bound on $V_2(H)$ is provided by
$V_1(H)$, and the lower bound by the one-step lookahead policy that takes exactly $T$ samples.}
\label{constrainedvalues}
\end{figure}

Recall in Section \ref{upperboundsec} we defined $V_1(H)$ as the optimal value of a 
policy with expected budget constraints, and $V_2(H)$ as the optimal value of a policy
with almost sure budget constraints. Also recall we proved $V_1(H)$ is equal to the solution of
the convex optimization program (\ref{upperboundeq}). In Figure \ref{expectedconstrainedvalue} 
we plot $V_1(H)$ as a function of the expected number of samples $T$, when $Y$ is a Brownian 
motion, with the same parameters as above. 
As discussed in Section \ref{upperboundsec}, $V_1(H)\ge V_2(H)$. 
A simple lower bound for $V_2(H)$ is the one-step lookahead policy
that takes exactly $T$ samples. For each $T\in\{1,\dots,10\}$ we estimated this lower bound
by simulating the one-step lookahead policy $50000$ times. In Figure \ref{tstepvalues} we
plot a region containing $V_2(H)$: the lower bound is provided by the one-step lookahead policy,
and the upper bound is provided by $V_1(H)$. The fact that shaded region in Figure \ref{tstepvalues}
is small means we have characterized $V_2(H)$ to high accuracy.
The fact that the lower bound provided by the one-step lookahead policy characterizes
$V_2(H)$ to such high accuracy shows that the one-step lookahead policy is very effective in the
constrained budget setting.

\section{Conclusion}
\label{conclusionsec}
In this paper, we consider a class of Bayesian optimization problems where
the underlying prior is a Markov process and we pay a cost for each sample. 
We show that the Bayes-optimal policy is computationally tractable, by way
of showing that the value function is completely determined by its values
on a $3$- or $4$-dimensional set. We use this optimal cost-per-sample policy to
compute the optimal value when there is no cost to sample, but there is a constraint
on the expected number of samples taken, as the result of a simple convex optimization
problem. We also use the optimal cost-per-sample policy to provide tight bounds
when the constraint on the number of samples taken is almost sure.
Computational experiments show that the optimal policy outperforms the
commonly used one-step lookahead policy, but also that the optimality gap between
one-step lookahead and the optimal policy is small, justifying the use of one-step
lookahead in practice. 

\appendix
\section{Proofs}
\label{proofsec}

\begin{proof}[\textbf{Proof of Lemma \ref{valueoutsideinterval}}]

We first show $V_{[a,b]}(H) \ge V_{[a,b]}(H^I)$.
Let $\sigma \in \Pi$.
Let $H^O=H\setminus H^I$ denote the set of initial observations outside $[a,b]$.
Define $\pi\in\Pi$ by
$\pi(K) \coloneqq \sigma(K\setminus H^O)$ for all $K\in\calH$.

Consider the following Markov processes:
\begin{enumerate} 
\item $(H_t^\pi)_{t\ge 0}$ with initial state $H_0=H$ operated under $\bbP^\pi$.
\item $(H_t^{I,\sigma})_{t\ge 0}$ with initial state $H_0^I=H^I$ operated under $\bbP^\sigma$.
\end{enumerate}
Now, note that 
\begin{equation}
R_{[a,b]}(H_\tau^\pi) = R_{[a,b]}(H_\tau^\pi\setminus H^O) \approx R_{[a,b]}(H_\tau^{I,\sigma})
\end{equation}
where the first equality holds (almost surely, under $\bbP^\pi$) because 
$Y$ is a Markov process and $H$ contains endoint observations,
and the second equality (in distribution) holds because
$H_\tau^\pi \setminus H^O$ under $\bbP^\pi$ is equal in distribution to
$H_\tau^{I,\sigma}$ under $\bbP^\sigma$. Moreover, $\tau^\pi$ under $\bbP^\pi$ starting from
initial state $H$ is equal in distribution to $\tau^\sigma$ under $\bbP^\sigma$ starting
from initial state $H^I$.
Thus
$\bbE^{\pi} \left[ R_{[a,b]}(H_{\tau}^\pi) - c\tau \mid H \right] =
\bbE^{\sigma} \left[ R_{[a,b]}(H_{\tau}^{I,\sigma}) - c\tau \mid H^I \right]$,
and the inequality $V_{[a,b]}(H) \ge V_{[a,b]}(H^I)$ is established.

To establish the reverse inequality, we apply the same logic as above. 
Let $\pi\in\Pi$ and define $\sigma\in\Pi$ 
by $\sigma(K) = \pi(K\cup H^O)$. Then the Markov processes $(H_t^\pi)_{t\ge 0}$
and $(H_t^{I,\sigma})_{t\ge 0}$ with initial states $H$ and $H^I$ respectively,
and operated under $\bbP^\pi$ and $\bbP^\sigma$ respectively,
share the same properties as before, but since $\sigma$ is now constructed
based on $\pi$ we conclude $V_{[a,b]}(H^I) \ge V_{[a,b]}(H)$.
\end{proof}

\begin{proof}[\textbf{Proof of Lemma \ref{nicepolicieslemma}}]
Fix $H_0\in \calH$.
Since $\bar{\Pi}_{[a,b]}\subseteq\Pi$ it follows that $\bar{V}_{[a,b]}(H_0)\le V_{[a,b]}(H_0)$.
We establish the three properties sequentially:
\begin{enumerate}
\item Suppose $\pi\in\Pi$ is such that $\bbP^\pi(\tau=\infty|H_0)\ne 0$. Then $\bbE^\pi[\tau|H_0]=\infty$.
Since $c>0$ it follows that $\per(\pi,c,H_0)=-\infty$ and so $\pi$ cannot be optimal. Thus the supremum
can be taken over $\Pi^1$.
\item Let $\pi\in\Pi$ be such that there is nonzero probability that, for some $t\ge 0$,
$\pi(H_t)\in H_t$. Write $x=\pi(H_t)$. Let $x_1,\dots,x_n$ be the complete set
of points that $\pi$ chooses to sample after sampling $x$ (note that all of these
points are random due to the randomness in the sample, except for $x_1$). 
Define $\sigma$ to be the same as $\pi$, except on $H_t$ where $\sigma$ samples $x_1$ first
and the distribution on the rest of the points is the same. Then 
$E^\pi[R_{[a,b]}(H_t)]=E^\sigma[R_{[a,b]}(H_t)]$ since sampling at $x$ again has no affect
on the reward, but $\sigma$ takes one fewer sample than $\pi$ and so it has better performance.
Repeating this process for every such $H_t$, one can construct a policy that never samples
a point that has already been sampled and that has better performance
than $\pi$. It follows that the supremum can be taken over $\Pi^1\cap\Pi^3$.
\item 
Now, let $\pi\in\Pi^1\cap\Pi^3$ and let $\mathcal{J}=\{H\in\calH : \pi(H)\notin [a,b]\}$.
Suppose $\mathcal{J}\ne\emptyset$ (equivalently, $\pi\notin\Pi^2_{[a,b]}$).
For each $H\in\calH$, let $x_1,\dots,x_{n_H}$ denote the random sequence of points $\pi$
samples until it chooses to stop sampling, or it samples inside $[a,b]$. That is,
if $\pi(H)\in [a,b]$ then $n_H$=1. If $\pi$ never samples inside $[a,b]$ after $H$,
then $x_{n_H}=\Delta$.
Define $\sigma(H)=x_{n_H}$.
It follows that $H_\tau^\sigma$ with initial state $H_0$ is equal in distribution
to $H_\tau^\pi\cap [a,b]$ also with initial state $H_0$. Thus
$\bbE^\sigma[R_{[a,b]}(H_\tau)\mid H_0]=\bbE^\pi[R_{[a,b]}(H_\tau)\mid H_0]$.
However, $\bbE^\sigma[\tau\mid H_0] \le \bbE^\pi[\tau \mid H_0]$.
Thus, the performance of $\sigma$ is equal or greater to the performance of $\pi$.
Hence the supremum can be taken over $\Pi^1\cap\Pi^2_{[a,b]}\cap\Pi^3$
and the result is established.
\end{enumerate}
\end{proof}

\begin{proof}[\textbf{Proof of Proposition \ref{translationthm}}]
By Lemma \ref{valueoutsideinterval} we assume all observations in $H$
are contained in $[a,b]$, i.e. $x\in [a,b]$ for all $(x,y)\in H$.

We show that for any policy $\pi\in\bar{\Pi}_{[a,b]}$ 
there exists a policy $\sigma\in\bar{\Pi}_{[a,b]}$ on $[a',b']$ such that
$\mathbb{E}^{\pi}\left[R_{[a,b]} - c\tau|H\right] = 
\mathbb{E}^{\sigma}\left[R_{[a',b']} - c\tau | H'\right].$

Fix any policy $\pi\in\bar{\Pi}_{[a,b]}$. Define $\sigma\in\bar{\Pi}_{[a',b']}$ by
$\sigma(K) \coloneqq T_{\ell} \circ \pi \circ T_{-\ell} (K\cap [a',b'])$ for
every $K\in\calH$. We use the intersection $K\cap [a', b']$ so that all observations
live at or above $0$, i.e. $T_{-\ell}(K\cap [a',b']) \in \calH$.

Now, consider the two Markov processes:
\begin{enumerate}
\item $\left(H_{t}\right)_{t\ge 0}$ under $\pi$ 
			with initial state $H_0=H$.
\item $\left(T_{-\ell}\left(H'_{t}\right)\right)_{t\ge 0}$ under $\sigma$ 
			with initial state $T_{-\ell}(H'_0) = T_{-\ell}(H')$.
\end{enumerate}
In (2), we apply the shift operator $T_{-\ell}$ to $H'_t$ so that
the two Markov Processes have the same initial state. 

We now show the two Markov processes have the same transition kernel.
Suppose $T_{-\ell}(H'_t)=H_t$ for some $t\ge 0$, that is, both
Markov Processes are in the same state at time $t$.  Note that $\pi$
and $T_{-\ell}\circ\sigma$ choose to sample the same point:
\[
T_{-\ell} \circ \sigma(H'_t) = 
T_{-\ell} \circ T_{\ell}\circ \pi \circ T_{-\ell}(T_{\ell}(H_t)\cap [a',b']) = 
 \pi (H_t \cap [a,b]) = \pi(H_t)
\]
where the final equality holds because
$H$ has all observations contained in $[a,b]$ and $\pi\in\bar{\Pi}_{[a,b]}$, so
$H_t$ must be contained in $[a,b]$ for all $t$.
Call this point $x_{t+1}$. It follows that every state $K$ with nonzero
probability for the $t+1$th time of both Markov Processes is of
the form $K=H_t \cup \{(x_{t+1},y)\}$ for some $y\in\bbR$. Then,
\begin{eqnarray}
\mathbb{P}^{\pi}\left(H_{t+1} = K | H_t\right) & = &
\mathbb{P}\left(Y(x_{t+1})\in dy|H_t\right) \\
& = & \mathbb{P} \left( Y(x_{t+1}) \in dy | T_{-\ell}(H'_t) \right) \\
& = &
\mathbb{P}^{\sigma} \left(T_{-\ell}(H'_t)=K|T_{-\ell}(H'_t)   \right).
\end{eqnarray}
Hence the two Markov Processes have the same transition kernel, and since
they have the same initial state, it follows they have the same distribution.

A simple consequence of this is that $\tau$ under $\pi$
and $\tau$ under $\sigma$ are identically distributed.
Indeed, $\tau_{\pi} \sim |H_{\tau}| - |H_0| \sim |H'_{\tau}| - |H'_0| \sim \tau_{\sigma}$.

Finally, observe that the reward is translation invariant: 
$R_{[a,b]}(K) = R_{[a',b']}(T_\ell(K))$ for any $K\in\calH$.
Thus,
\begin{eqnarray}
\label{reward_firstline}
\bbE^{\pi}\left[ R_{[a,b]}(H_\tau) \mid H_0\right] & = &
\bbE^\pi \left[ R_{[a',b']}(T_\ell(H_\tau) \mid H_0\right] \\
\label{reward_secondline}
 & = & 
\bbE^\pi \left[ R_{[a',b']}(T_\ell(H_\tau)) \mid T_\ell(H_0) \right] \\
\label{reward_thirdline}
& = &
\bbE^\sigma \left[ R_{[a',b']}(H'_\tau) \mid H'_0\right] 
\end{eqnarray}
where the equality between (\ref{reward_firstline}) and (\ref{reward_secondline})
holds because $T_\ell$ is a bijection, and equality between
(\ref{reward_secondline}) and (\ref{reward_thirdline}) holds because
$T_\ell(H_\tau)\mid T_\ell(H_0)$ under $\pi$ is equal in distribution
to $H'_\tau\mid H'_0$ under $\sigma$.

Thus $\mathbb{E}^{\pi}\left[R_{[a,b]} - c\tau|H\right] = 
\mathbb{E}^{\sigma}\left[R_{[a',b']} - c\tau | H'\right]$, as we
set out to show. It follows that
$V_{[a,b]}(H)\le V_{[a',b']}(H')$.
Setting $a\coloneqq a+\ell$, $b\coloneqq b+\ell$ and $\ell \coloneqq -\ell$ establishes
the reverse inequality, and equality follows.
\end{proof}

\begin{doublespace}
\end{doublespace}

\bibliography{writeup}
\bibliographystyle{plain}

\end{document}